\documentclass[12pt]{article}

\usepackage{ifthen}
\newcommand{\plain}{1} %Plain

\ifthenelse{\plain = 1}
{}
{\smartqed}

\usepackage{graphicx}
\usepackage{eucal,amssymb,amsmath,graphicx, amsfonts, latexsym}

\newcommand{\bs}{\boldsymbol}

\newcommand{\Barabasi}{Barab{\'a}si}
\newcommand{\Erdos}{Erd\H{o}s}
\newcommand{\Lageras}{Lager{\aa}s}
\newcommand{\Renyi}{R\'{e}nyi}
\newcommand{\Poi}{\mathrm{Poi}}
\newcommand{\pmaj}{p_{\mathrm{maj}}}

\ifthenelse{\plain = 1}
{
\setlength{\textheight}{24cm}\setlength{\textwidth}{16cm}
\setlength{\topmargin}{-2cm} \setlength{\oddsidemargin}{-0.2cm}
\setlength{\evensidemargin}{0.0cm}
\newcommand{\hfigwidth}{7cm}
\newcommand{\figwidth}{10cm}
}
{
\newcommand{\hfigwidth}{0.45\textwidth}
\newcommand{\figwidth}{0.75\textwidth}
}

\begin{document}

\ifthenelse{\plain=1}
%Plain title stuff
{
\parskip=5pt plus1pt minus1pt \parindent=0pt
\title{A network with tunable clustering, degree correlation and degree distribution, and an epidemic thereon}
\author{Frank Ball\thanks{School of Mathematical Sciences, University of Nottingham, University Park,
Nottingham NG7 2RD, UK. Email: {\tt frank.ball@nottingham.ac.uk}}, Tom Britton\thanks{Department of Mathematics, Stockholm University, 106 91 Stockholm, Sweden. Email: {\tt tomb@math.su.se}} and David Sirl\thanks{Mathematics Education Centre, Loughborough University, Loughborough, Leicestershire, LE11 3TU, UK. E-mail: {\tt d.sirl@lboro.ac.uk}}}
\date{6th July 2012}
}
%JMB title stuff
{
\title{A network with tunable clustering, degree correlation and degree distribution, and an epidemic thereon
%\thanks{Grants or other notes
%about the article that should go on the front page should be
%placed here. General acknowledgments should be placed at the end of the article.}
}
%\subtitle{Do you have a subtitle?\\ If so, write it here}

%\titlerunning{Short form of title}        % if too long for running head
\author{Frank Ball \and Tom Britton \and David Sirl} %\authorrunning{Ball \etal}
\institute{F.~G.~Ball \at School of Mathematical Sciences, University of Nottingham, University Park, Nottingham NG7 2RD, UK. \email{frank.ball@nottingham.ac.uk} \and
T.~Britton \at Department of Mathematics, Stockholm University, SE-106 91 Stockholm, Sweden. \email{tomb@math.su.se} \and
D.~J.~Sirl \at Mathematics Education Centre, Loughborough University, Loughborough, Leicestershire LE11 3TU, UK. \email{d.sirl@lboro.ac.uk}}
%\authorrunning{Ball, Britton and Sirl} % if too long for running head
\date{Received: date / Accepted: date}
% The correct dates will be entered by the editor
}
\maketitle

\begin{abstract}
A random network model which allows for tunable, quite general forms of clustering, degree correlation and degree distribution is defined. The model is an extension of the configuration model, in which stubs (half-edges) are paired to form a network. Clustering is obtained by forming small completely connected subgroups, and positive (negative) degree correlation is obtained by connecting a fraction of the stubs with stubs of similar (dissimilar) degree. An SIR (Susceptible $\to$ Infective $\to$ Recovered) epidemic model is defined on this network. Asymptotic properties of both the network and the epidemic, as the population size tends to infinity, are derived: the degree distribution, degree correlation and clustering coefficient, as well as a reproduction number $R_*$, the probability of a major outbreak and the relative size of such an outbreak. The theory is illustrated by Monte Carlo simulations and numerical examples.  The main findings are that clustering tends to decrease the spread of disease, the effect of degree correlation is appreciably greater when the disease is close to threshold than when it is well above threshold and disease spread broadly increases with degree correlation $\rho$ when $R_*$ is just above its threshold value of one and decreases with $\rho$ when $R_*$ is well above one.
\ifthenelse{\plain=1}{}
%JMB keywords, MSC
{
% \PACS{PACS code1 \and PACS code2 \and more}
}
\end{abstract}

\ifthenelse{\plain=1}
%Plain Keywords, MSC
{

}
{}

\section{Introduction}

Ever since the pioneering work of \Erdos\ and \Renyi\ (1959) on a simple random graph there have been numerous important contributions on random graph models with the aim of making them more flexible and realistic. For example, the configuration model (Molloy and Reed (1995) and Newman et al.~(2001))  defines a network allowing for more or less arbitrary \emph{degree distribution} $F_D$, the distribution describing the number of neighbours $D$ of a randomly selected node (which in the epidemic context represents an individual) in the network. (For simplicity, from now on we refer to $D$ as the degree distribution.)  This extension was important for two reasons: most empirical networks tend to have much heavier tailed  degree distributions than the Poisson distribution of the \Erdos-\Renyi\ (E-R) network, and networks with heavy tail degree distributions have been shown to exhibit rather different properties when compared with the E-R network; for example, if an epidemic outbreak takes place on the network the epidemic threshold $R_0$ is much higher (or even infinite) as compared to the same epidemic taking place on an E-R network with the same mean degree (Andersson  (1999)).

Two other properties of real world networks that are not present in E-R networks are \emph{clustering} and \emph{degree correlation}. The clustering coefficient $c$ measures how likely it is that two neighbours of a randomly selected node are neighbours themselves. The E-R network has no clustering whereas nearly all empirical networks have positive clustering, with typical values in the range 0.1--0.5 out of the possible range 0--1 (see Newman (2003), Table 3.1).  The degree correlation $\rho$ instead measures the correlation between the degrees of the adjacent individuals of a randomly selected \emph{edge}. The E-R network has $\rho =0$ whereas `random' networks with heavy tail degree distribution tend to have $\rho >0$ (van der Hofstad and Litvak (2012)). Empirical networks, on the other hand, have both positive and negative degree correlation: there seems to be a tendency for computer networks to have $\rho <0$ whereas social networks (our main interest) typically have $\rho>0$ (see Newman (2003), Table 3.1). There are numerous network models studied in the literature, with the aim of allowing one or several of these three extensions (of \emph{local} properties) from the original E-R-network (see some references below where the focus is also on epidemics evolving on the network); the term `local' refers to the fact that it is sufficient to observe nodes and their neighbourhoods to determine/estimate such properties (the complete network need not be observed in order to evaluate them). The current paper defines a model in which $D$, $c$ and $\rho$ can be made more or less arbitrary.

There are of course other important extensions in addition to allowing for arbitrary degree distribution, degree correlation and clustering. Further local properties considered in many models for social networks are households and other fully connected smaller units (e.g.\ Ball et al.~(1997)), and models in which nodes and/or edges are of different types (e.g.\ Britton et al.\ (2007), Ball and Sirl (2012)). Several models have also been proposed which combine household and network structure, for example Trapman (2007), Gleeson (2009), Ball et al.\ (2010) and Ma et al.\ (2012). Other models aim to study and extend the range of global properties, such as small world networks (Watts and Strogatz (1998)) and dynamic network models (\Barabasi\ and Albert, (1999)). This paper does not address these (or any other) extensions; the focus being on degree distribution, degree correlation and clustering.

Our main motivation for studying networks is to investigate social networks and to examine what effect the three above-mentioned properties have in the event of an infectious disease entering the community; both in terms of the possibility and probability of an epidemic outbreak taking off, and also how large such an outbreak will be if it does take off. We study the class of SIR epidemics (e.g.~Andersson and Britton (2000)) in which individuals are at first Susceptible (except for some introductory infectious cases) and those who get infected become Infectious for a random period of time when they may infect their network neighbours, after which they Recover and become immune to further infection. See, for instance, Diekmann et al.~(1998), Andersson (1999) and Diekmann and Heesterbeek (2000, Ch.\ 10) for early analytical contributions in this area.

As mentioned above there have been many contributions to this area of research, in particular over the last decade or two. Allowing for arbitrary degree distribution, and studying its effect on an epidemic, dates back longer. May and Anderson (1987) concluded (when modelling the spread of HIV) that a heavy tail degree distribution makes the reproduction number $R_0$ large or even infinite. The important insight from their analysis was that diseases with very low transmission probability still may be at risk of epidemics taking off in networks having small mean degree, if the \emph{variance} of the degree distribution is very large. The effect of clustering on epidemics has been studied in, for example, Britton et al.\ (2008), Miller (2009) and Newman (2009). Degree correlation has often been analysed in combination with clustered networks (e.g.\ Gleeson et al.~(2010)). The impact of clustering and degree correlation on epidemics on networks has been studied empirically using simulation by Badham and Stocker (2010) and Isham et al.~(2011). The main focus of most papers concerning epidemics on networks with controllable clustering, degree correlation and/or degree distribution lies in studying how these features affect the basic reproduction number $R_0$, i.e.~the \emph{possibility} of having an major outbreak. To derive the \emph{probability} of such an outbreak, and its likely size in the event that it takes off, requires significantly deeper analysis; which for several of the above-mentioned models still is missing.

The current paper introduces a network model which (i) allows for more or less arbitrary clustering, degree correlation and degree distribution, and (ii) permits theoretical analysis of epidemics defined on the network.  As in the configuration model, the network is formed by attaching stubs (i.e.~half-edges) to individuals, which are then paired to form the edges of the network. The degree of an individual is the number of stubs emanating from it.  The desired clustering and degree distribution is obtained by having two types of stubs going out from individuals. A fraction of stubs is local (which fraction being closely related to the desired clustering); the remaining stubs are global and are connected randomly (as described below) among stubs from all individuals. The local stubs are connected by grouping individuals into small local groups (`households'). For example, an individual with four local stubs is connected to four other individuals having local degree 4, thus forming a group of 5 completely connected individuals (contributing to increased clustering). The degree distribution is given by the distribution of the sum of the local and global degree of a typical individual. Finally, the desired degree correlation $\rho$ is obtained by manipulating how the global stubs are connected, which is controlled by a parameter $r$ satisfying $-1 \le r \le 1$.  With probability $1-|r|$ a stub is connected uniformly at random among all global stubs. With probability $|r|$ the stub is connected to a stub having very similar total degree (if $r>0$) or `opposite' total degree (if $r<0$). 

The remainder of the paper is organised as follows.  A more rigorous definition of the model appears in Section~\ref{sec-net-epi}, where a continuous-time SIR epidemic on the network is also defined. In Section~\ref{sec-net-prop}, we derive expressions for the degree distribution $D$, clustering coefficient $c$ and degree correlation $\rho$, as  functions of the model parameters, and discuss the more relevant reverse problem of choosing model parameters to obtain a desired $c$, $\rho$ and $D$, using a Poisson total degree distribution as a template.  We also describe a simple rewiring algorithm, motivated by Miller (2009) and Gleeson et al.~(2010), which permits the clustering in a network to be reduced in a controlled fashion without changing $\rho$ or $D$.  In Section~\ref{epid-prop}, we analyse the main characteristics of epidemics defined on the network for suitably large population sizes, by exploiting approximating branching processes.   Specifically, in Section~\ref{estab}, we obtain a threshold parameter $R_*$ which determines whether or not a major outbreak is possible, and derive the probability that a major outbreak occurs (assuming that the infectious period is constant) and, in Section~\ref{zglobal}, we derive the relative final size (i.e.~the proportion of the population that is ultimately removed) of a major outbreak.  In Section~\ref{rewired}, we describe how these results on epidemics are modified to incorporate rewiring and prove that, if all other parameters are held fixed, such rewiring increases the threshold parameter $R_*$ and both the probability and relative final size of a major outbreak.  In Section~\ref{numexamples} we illustrate the theory with some numerical examples which demonstrate that the effect of degree correlation on epidemic properties is appreciably greater when the disease is just above threshold than when it is well above threshold.  Moreover, both the probability and size of a major outbreak broadly increase with $\rho$ when the disease is just above threshold, while they broadly decrease with $\rho$ when the disease is well above threshold.  However, this behaviour is not monotonic, particularly when clustering is low and $R_*$ is close to one.  We conclude with a brief discussion in Section~\ref{sec-disc}.

\section{The network model and the epidemic}
\label{sec-net-epi}

\subsection{The network model}\label{sec-net-mod}

Consider a network of undirected edges with $n$ nodes (individuals). Below we define how to construct the network. First we define a set of random variables and briefly explain their interpretation in the network.

Let $G$ be a discrete non-negative random variable with distribution $\{p_k\}$ referred to as the `global degree', let $H$ be another strictly positive discrete random variable with distribution $\{\pi_h\}$. In some cases $H$ will reflect the household distribution in the community, but in applications where the underlying network has no household structure $H$ is simply a device to introduce clustering into the network.
%will only have probability mass at 1 and 3: $P(H=3)=\pi_3= \alpha=1-\pi_1=1-P(H=1)$, and $\alpha$ denotes the fraction of ''local triangles'' as compared to the ''local singletons''. Note that this implies that a node belongs to a local triangle (rather than being a local singleton) with probability $3\alpha/(1-\alpha + 3\alpha)=3\alpha/(1+2\alpha)$. 
Finally, let $r$ be a real number satisfying $-1\le r\le1$. The value of $|r|$  reflects how often outgoing global edges connect to nodes of similar (if $r>0$) or `opposite' (if $r<0$) `total degree'. Let $X$ be a Bernoulli random variable with parameter $|r|$, so $P(X=1)=|r|=1-P(X=0)$, this variable will determine if a stub will connect to a random stub or a stub with similar/`opposite' degree.

The network is constructed as follows. Let $H_1, H_2, \cdots $ be independent and identically distributed copies of the random variable $H$. Label the $n$ nodes $1,2,\cdots,n$ and group the first $H_1$ nodes into local group (household) one, nodes $H_1+1, H_1+2, \cdots, H_1+H_2$ into group 2 and so on until all individuals belong to a local group (the last group will have a `truncated' size). All nodes of a local group are connected to each other (for example, the first $H_1$ nodes make up a fully connected component with all individuals having local degree $H_1-1$). Let $G_1, G_2,\cdots ,G_n$ be independent and identically distributed copies $G$; $G_i$ denotes the global degree of node $i$. The total degree of individual $i$, $D_i$, equals the global degree plus the local degree, the local degree being one less than the group size. For example, a node residing in a local triangle ($H=3$) and having global degree $k$ has total degree $k+2$, whereas a local singleton ($H=1$) with global degree $j$ has total degree $j$. A node having global degree $k$ has $k$ outgoing stubs, and each of these stubs is labelled with an independent copy of $X$ (stubs having independent and identically distributed $X$-variables with $P(X=1)=|r|=1-P(X=0)$) and the \emph{total} degree of the node from which it emanates. All outgoing stubs in the network with label $X=0$ are connected pairwise completely at random. The remaining stubs (having $X=1$) are also connected randomly but in a different manner. This is done by ordering all global stubs having label $X=1$ (suppose that there are $n_1$ such stubs) according to their total degree, and then separating the empirical distribution of global degrees so generated into $n_Q$ (a fixed and freely chosen positive integer) equally sized quantiles. (If $n_1/n_Q$ is not an integer then the $n_Q$ quantiles are made as equal in size as possible.)  The first such quantile hence consists of the $n_1/n_Q$ stubs having smallest label (i.e.\ total degree) and so on.  If $r>0$, each quantile is treated in turn and all the stubs in that quantile are paired uniformly at random.  If $r<0$, the stubs in the first quantile are paired uniformly at random with those in the $n_Q$th quantile, the stubs in the second quantile are paired uniformly at random with those in the $n_Q-1$th quantile, and so on.  Thus, if $n_Q$ is odd, the stubs in the middle quantile are paired uniformly at random with each other. The effect of this pairwise connection is that nodes of similar total degree will be connected if $r>0$, whereas nodes of rather different total degree will be connected if $r<0$; in both cases leading to correlated degrees (but of different sign).  There may be one unattached stub having label $X=0$ and at most $n_Q$ unattached stubs having label $X=1$ following the above pairings.  These are simply ignored.  This has no effect on the asymptotic properties of the network, nor on epidemics defined thereon, as $n \to \infty$.  In the above construction, all the $H,G$ and $X$ random variables are assumed to be independent.

The network is hence made up of local completely connected groups having groups size distribution $\{ \pi_h\}$ (as $n$ goes to infinity the effect of the last group having a truncated household size is negligible). On top of this, each individual has global edges, the number being distributed as $G$. Some of these will be formed by connecting to other random stubs, the others will be formed by connecting to other stubs having similar or `opposite' degree, thus creating positive or negative degree correlation. The construction of global edges may result in the presence of multiple edges and self-loops. However, if the degree distribution $D$ has finite variance, the fraction of these will be negligible as $n\to\infty$, so removing them has negligible effect on the degree distribution and how stubs are connected (cf.\ Durrett (2006, Theorem 3.1.2) and Janson (2009)). The special case where $r=0$ or $n_Q=1$ is the network and households model (without degree correlation beyond that induced by the presence of households) studied by Ball et al.~(2010), since in either of these situations all global stubs are simply paired uniformly at random.

\subsection{An epidemic model on the network}\label{sec-epid}

We now define a continuous-time epidemic model for the spread of an SIR-type infectious disease upon the network defined in Section~\ref{sec-net-mod}. We suppose that there is one initial infective, chosen uniformly at random from the $n$ individuals (nodes) in the population and that the remainder of the population is susceptible.  The infectious periods of different infectives are each distributed according to a random variable $I$, having an arbitrary but specified distribution.  Throughout its infectious period, a given infective makes infectious contacts with any given neighbour (either local or global) in the network at the points of a homogeneous Poisson process having rate $\lambda$.  A susceptible becomes infective as soon as it is contacted by an infective and an infective becomes removed (and plays no further part in the epidemic) at the end of its infectious period.  Contacts between an infective and an infective or removed individual have no effect.  All Poisson processes describing infectious contacts (whether or not either or both individuals involved are the same) and all infectious periods are mutually independent; they are also independent of the random variables used to construct the network.  The epidemic ends when there is no infective remaining in the population.

\section{Properties of the network model}\label{sec-net-prop}

We now derive the total degree distribution $D$, the clustering coefficient $c$ and the degree correlation $\rho$ for the network defined in Section~\ref{sec-net-mod}. We treat the asymptotic case where the number of nodes $n$ tends to infinity.

\subsection{The degree distribution}

We start with the degree distribution. From the construction it follows immediately that a node has global degree $G$. The local degree is one less than the household size, and the household size of a randomly selected node has distribution $\{\tilde \pi_h\}$, where $\tilde \pi_h=h\pi_h/\mu_H$ and $\mu_H=\sum_jj\pi_j$, i.e.\ the size-biased local group-size distribution. Let $\tilde H$ denote a random variable having the size-biased household distribution. It then follows that the total degree distribution (in the network) is given by
\begin{equation}
D \overset{D}{=}G+\tilde H-1,\label{deg-dist}
\end{equation}
where $\overset{D}{=}$ means equal in distribution and $G$ and $\tilde{H}$ are independent. In particular it follows that the mean total degree is
\[
 \mu_D=\mu_G+\frac{\sigma_H^2}{\mu_H}+\mu_H - 1.
\]
(Throughout the paper, for a random variable, $X$ say, $\mu_X$ and $\sigma_X^2$ denote respectively the mean and variance of $X$.)

\subsection{The clustering coefficient}\label{sec-clust}

There are several measures of clustering used in the literature.   We use a `probabilistic' one (see, for example, Trapman (2007)) where an ordered triplet of nodes $(i,j,k)$ is selected completely at random among all such ordered triplets for which $i$ is directly connected to $j$ and $j$ is directly connected to $k$. The clustering coefficient $c$ is then defined as the probability that $i$ and $k$ are also directly connected (i.e.\ that $i$, $j$ and $k$ form a triangle). Thus $c$ is given by the fraction of ordered triplets in the network that are triangles. The clustering coefficient of the present network model is identical to that of the model in Ball et al.~(2010), since the models differ only in the way that global stubs are paired.  For large $n$, the proportion of ordered triangles that are not wholly within households is small and zero in the limit as $n\to \infty$. Thus, asymptotically, the global pairings do not yield triangles in either of the two models, explaining why the clustering coefficients are the same for the two models. Hence, from equation (14) of Ball et al.~(2010), the clustering coefficient $c=c(G,H,r)$ is given by
\begin{equation}
c=\frac{{\rm E}[H(H-1)(H-2)]}{{\rm E}[(H(G+H-1)(G+H-2)]},\label{clust}
\end{equation}
where $G$ and $H$ are the household and global degree distributions of the network.

\subsection{The degree correlation}\label{sec-corr}

We now formulate an expression for the degree correlation $\rho$ of the current network model. One way to define $\rho$ is to pick a random edge in the network and let $\rho$ be the correlation between the total degrees of the nodes adjacent to this edge (Newman, 2002a).  The derivation of $\rho$ involves long but standard computations which are given in the appendix. A key step in the derivation is to first condition on whether the chosen edge is a global or a local edge, the former having probability $p_G$ given by
\begin{equation}
p_G= \frac{\mu_G}{\mu_G+\mu_{\tilde H}-1}.\label{frac-glob}
\end{equation}
If the edge is global the degree covariance (of the right and left node adjacent to the edge) comes from the two stubs having the same (or `opposite') quantile(s), which happens with probability $|r|$, and if the edge is local the degree covariance stems from the nodes having the same local degree.

Before giving the expression for the degree correlation $\rho=\rho(G,H,r)$ some more notation is required.  Let $\hat H$ denote a random variable giving the household size of a household edge chosen uniformly at random from all household edges.  Since a household of size $h$ contains $\binom{h}{2}$ edges, ${\rm P}(\hat H=h) \propto \binom{h}{2} \pi_h$ ($h=2,3,\cdots$), so
\begin{equation*}
{\rm P}(\hat H=h)=\frac{h(h-1)\pi_h}{{\rm E}[H(H-1)]} \quad (h=2,3,\cdots).
\end{equation*}

Let $\tilde{D}$ and $\tilde{Q}$ denote respectively the \emph{total} degree and quantile of a stub chosen uniformly at random from all stubs in the limit as $n \to \infty$.  Then $\tilde{D}\overset{D}{=}\tilde{G}+\tilde{H}-1$, where $\tilde{G}$ and $\tilde{H}$ are independent, and
$\tilde{G}$ denotes a random variable having the size-biased global degree distribution $\{\tilde{p}_g\}$, where $\tilde{p}_g=g p_g/\mu_G$ $(g=1,2,\cdots)$.
For $i=1,2,\cdots,n_Q$ and $d=1,2,\cdots$, let $p_{\tilde{Q}|\tilde{D}}(i|d)={\rm P}(\tilde{Q}=i|\tilde{D}=d)$ and $p_{\tilde{D}|\tilde{Q}}(d|i)={\rm P}(\tilde{D}=d|\tilde{Q}=i)$.  (These conditional probabilities are derived easily from the probability mass function of $\tilde{D}$, noting that if $\tilde{u}_0=0$ and $\tilde{u}_d={\rm P}(\tilde{D}\le d)$ $(d=1,2,\cdots)$ then ${\rm P}(\tilde{D}=d,\tilde{Q}=i)=\max\left\{\min(\tilde{u}_d,\frac{i}{n_Q})-\max(\tilde{u}_{d-1},\frac{i-1}{n_Q}),0\right\}$ $(d=1,2,\cdots; i=1,2,\cdots,n_Q)$.)  Define the function $g_{\tilde{D},n_Q}(r)$ by
%\marginal{Is there an interpretation for $g_{\tilde{D},n_Q}(r)$?} Yes but complicated: g=|r| g_1, where g_1 is the covariance of conditional expectations in equation \eqref{covEXLEXRIGQ}
\begin{equation}
\label{gtilded}
g_{\tilde{D},n_Q}(r)=\begin{cases}
r\left(\frac{1}{n_Q}\sum_{i=1}^{n_Q}(\mu_{\tilde{D}}^{(i)})^2-\mu_{\tilde{D}}^2\right) & \mbox{if } r \ge 0, \\
|r|\left(\frac{1}{n_Q}\sum_{i=1}^{n_Q}\mu_{\tilde{D}}^{(i)}\mu_{\tilde{D}}^{(n_Q+1-i)}-\mu_{\tilde{D}}^2\right)& \mbox{if } r<0,
\end{cases}
\end{equation}
where
\begin{equation}
\label{mudtildei} 
\mu_{\tilde{D}}^{(i)}=\sum_{d=1}^{\infty}d p_{\tilde{D}|\tilde{Q}}(d|i) \quad
(i=1,2,\cdots,n_Q).
\end{equation}
It is shown in the appendix that
\begin{equation}
\label{degcorr}
\rho=\frac{(1-p_G)\sigma_{\hat{H}}^2+p_G g_{\tilde{D},n_Q}(r)+p_G(1-p_G)\left(\mu_{\hat{H}}-\mu_{\tilde{H}}-\frac{\sigma_G^2}{\mu_G}\right)^2}
{(1-p_G)\left(\sigma_{\hat{H}}^2+\sigma_G^2 \right)+p_G\left(\sigma_{\tilde{H}}^2+\sigma_{\tilde{G}}^2\right)
+p_G(1-p_G)\left(\mu_{\hat{H}}-\mu_{\tilde{H}}-\frac{\sigma_G^2}{\mu_G}\right)^2}.
\end{equation}

\subsection{Rewiring}
\label{rewiring}
Note that for household size and global degree distributions $H$ and $G$, the degree distribution $D$ and the clustering coefficient $c$ are both independent of the parameter $r$.  Thus, by letting $r$ vary between $-1$ and $+1$ and keeping the distributions of $H$ and $G$ fixed, it is straightforward to tune the degree correlation in our network model without changing the degree distribution or clustering coefficient of the network.  However, if we keep $r$ fixed and vary, for example, the household size distribution to change the clustering coefficient of the network, then its degree distribution $D$ and degree correlation $\rho$ change also.  This observation means that it is more difficult to tune just the clustering coefficient in a network.  One way around this problem is to extend the rewiring construction of Gleeson et al.~(2010) (see also Miller (2009), where the idea first originated) to our model.

Suppose that we construct a realisation of our network model and then colour all global edges green and all household edges red. Household edges are also labelled according to their household size.  Let $p_{RW}$ be a real number satisfying $0 \le p_{RW}\le1$.  Then, independently for each household, with probability $p_{RW}$ the red edges in a household are each broken into two stubs, which retain their colour and household-size labels. For each $h=2,3,\cdots$, the red stubs with label $h$ are now joined uniformly at random, which, together with the green edges and unbroken red edges creates a new network.

Observe that the above rewiring does not alter the degree distribution or the correlation structure (and in particular the degree correlation) of the network but it does change its clustering coefficient.  Let $c(G,H,r,p_{RW})$ denote the clustering coefficient for the model with rewiring probability $p_{RW}$, so $c(G,H,r,0)$ is the clustering coefficient of our model without rewiring.  In the limit as $n \to \infty$, the proportion of triangles that are not wholly within unbroken households tends to zero, whence $c(G,H,r,p_{RW})=(1-p_{RW})c(G,H,r,0)$.  Thus, given our network model without rewiring, it is straightforward to use the above rewiring to tune the clustering coefficient to be any value between $0$ and that of the model without rewiring.

\subsection{Tuning}
\label{Sec-Tuning}
The formulae given in Sections \ref{sec-clust} and \ref{sec-corr} are fairly long but simplify appreciably for the special situation where both the household sizes and the global degrees follow Poisson-based distributions.  Specifically, suppose that, with 
$0 \le \mu < \gamma$, $G$ follows a Poisson distribution with mean $\gamma-\mu$, which we denote by $\Poi(\gamma-\mu)$, and $H$ follows a Poisson distribution with mean $\mu$ that is conditioned on being strictly positive, which we denote by $\Poi^+(\mu)$. Here we interpret $\Poi^+(0)$ to be $\lim_{\mu\downarrow0}\Poi^+(\mu)$, the distribution identically equal to 1. Thus $\pi_h=(1-{\rm e}^{-\mu})^{-1}\mu^h{\rm e}^{-\mu}/h!$ $(h=1,2,\cdots)$.  Then $\tilde{H}-1\sim \Poi(\mu)$ and it follows from~\eqref{deg-dist} that the total degree $D\sim\Poi(\gamma)$.  Further,
$1-p_G=\mu/\gamma$ and $\hat{H}-2\sim\Poi(\mu)$, so using~\eqref{clust} and \eqref{degcorr}, the formulae for the clustering and degree correlation are given by:
\begin{equation}
\label{Poicrho}
c=\left(\frac{\mu}{\gamma}\right)^2\qquad \mbox{and}\qquad \rho=\frac{1}{\gamma^2}\left[\mu^2+(\gamma-\mu)g_{\gamma,n_Q}(r)\right],
\end{equation}
where $g_{\gamma,n_Q}(r)$ is given by~\eqref{gtilded} with $\tilde{D}\sim 1+\Poi(\gamma)$.

\begin{figure}
\begin{center}
\resizebox{\figwidth}{!}{\includegraphics{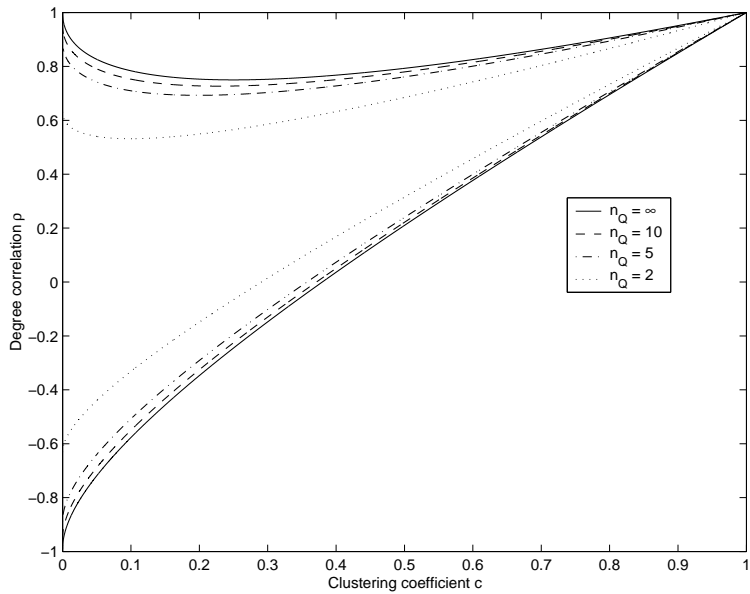}}
\end{center}
\caption{Plot showing bounds on possible values of $(c,\rho)$ when $D\sim\Poi(10)$.}
\label{fig:poissontune}
\end{figure}

Observe that $g_{\gamma,n_Q}(0)=0$, so $c=\rho$ when $r=0$, i.e.~for the model studied in Ball et al.~(2010), Sections 4.3 and 4.4.  Suppose that $\gamma$ and $\mu$ are held fixed, so the clustering coefficient $c$ is also held fixed.  Then as $r$ varies from $-1$ to $+1$ the degree correlation $\rho$ varies between the values obtained by setting $r=-1$ and $r=1$ in the formula for $\rho$ in~\eqref{Poicrho}. These lower and upper values for $\rho$ are shown in Figure~\ref{fig:poissontune} as functions of $c$ for different choices of the number of quantiles $n_Q$, for the case when $\gamma=10$.  In the limit as $n_Q\to\infty$, if $r>0$ then a stub with label $X=1$ is paired, almost surely, with a stub having the same total degree and $g_{\gamma,n_Q}(1) \to {\rm var}(\tilde{D})=\gamma$ (recall $\tilde{D}\sim 1+\Poi(\gamma)$).  It follows that the corresponding upper value for $\rho$ is $1+c-\sqrt{c}$.  In the same limiting situation, if $r<0$ then a stub with label $X=1$ is paired, almost surely, with a stub having the `opposite' total degree.  There is no simple expression for $\lim_{n_Q \to \infty}g_{\gamma,n_Q}(-1)$, though it is easily computed.  Observe from Figure~\ref{fig:poissontune} that very little extra is gained, in terms of the range of possible $(c,\rho)$, by choosing a large value of $n_Q$. In practice, a small value of $n_Q$ is beneficial as the proportions of self-loops and parallel edges between nodes, resulting from the pairing of stubs, both increase with $n_Q$. Additionally, large values of $n_Q$ mean that the approximating branching processes have many types and numerical calculation of quantities of interest becomes more computationally intensive.

Write $c=c(\gamma,\mu,r)$ and $\rho=\rho(\gamma,\mu,r)$ to show explicitly their dependence on the parameters and, for $\gamma>0$, let $A_{\gamma}=\{(c(\gamma,\mu,r),\rho(\gamma,\mu,r)):0 \le \mu \le \gamma, -1 \le r \le 1\}$ be the set of possible values $(c,\rho)$ in our model when the total degree is $\Poi(\gamma)$.  For any $(c,\rho)\in A_{\gamma}$ there is a unique $(\mu,r)$ such that $(c(\gamma,\mu,r),\rho(\gamma,\mu,r))=(c,\rho)$, so the model without rewiring can be tuned uniquely to any attainable $(c,\rho)$. 
%If we allow rewiring, it is easily seen that by choosing the rewiring probability $p_{RW}$ appropriately, for each $(c,\rho)$ that is in the interior of $A_{\gamma}$ there is a continuum of models with clustering coefficient $c$ and degree correlation $\rho$.   The same is also true for any $(0,\rho) \in A_{\gamma}$.
If we allow rewiring, it is easily seen that by choosing the rewiring probability $p_{RW}$ appropriately, for each $(c,\rho)$ lying strictly above the lower boundary of $A_{\gamma}$, there is a continuum of models with clustering coefficient $c$ and degree correlation $\rho$.

A similar analysis to the above holds for other choices of total degree distribution $D$, though note that not all distributions $D$ can be decomposed as in~\eqref{deg-dist} in such a way that the clustering may be tuned continuously.  Distributions $D$ for which this is possible include negative binomial and compound Poisson.  Indeed any distribution $D$ that is infinitely divisible may be decomposed so that the clustering coefficient is any rational number in $[0,1)$.

\section{Epidemics on network without rewiring}
\label{epid-prop}
\subsection{Establishment of the epidemic}
\label{estab}
\subsubsection{Approximating forward branching process}
\label{forwardBP}
The initial infective triggers a local (i.e~within-household) epidemic in its household.  Each infective in that local epidemic (including the initial infective) may make (global) infectious contact with individuals in other households.  If the population size $n$ is large, the probability that such global infectious contacts are all with individuals in previously uninfected households is close to one, owing to the random way in which the underlying network is formed.  It follows that in the early stages of an epidemic the process of infected households may be approximated by a branching process, with individuals in the branching process corresponding to infectious households in the epidemic process.  Unless $r=0$ or $n_Q=1$, this branching process needs to be multitype, since the degrees of endpoints of a global edge with $X=1$ are correlated.  Except for the ancestor, the type of an individual in the branching process is obtained by considering the primary infective, $i^*$ say, in the corresponding single-household epidemic.  The type of the individual is given by the total-degree quantile of the stub used in constructing the global edge along which $i^*$ was infected in the epidemic.  Thus there are $n_Q$ types of individual in the branching process.  The ancestor of the branching process is not typed in this fashion since the initial infective in the epidemic is chosen uniformly at random from the population and not infected along a global edge in the network.  Nevertheless, the offspring distribution of the ancestor in the branching process depends on the household size and global degree of the initial infective in the epidemic.

Following Ball et al.~(2009), the above branching process is termed a forward branching process as it approximates the forward spread of an epidemic process.  In Section~\ref{zglobal} we consider a backward branching process, which approximates an inverse epidemic process.

The approximation of the early stages of the epidemic process by the forward  branching process can be made precise by constructing the branching process and, for each $n=1,2,\cdots$, a realisation of the epidemic process on a common probability space and using a coupling argument to show that, as $n \to \infty$, the process of infected households in the epidemic process converges almost surely to the multitype branching process; cf.~Ball and Sirl (2012).  Thus, if the population size $n$ is sufficiently large, the probability that the epidemic becomes established and leads to a major outbreak is given approximately by the probability that the branching process survives (i.e.~does not go extinct).
Moreover, whether or not a major outbreak can occur with non-zero probability is determined by whether or not the branching process is supercritical.

We now determine the means and probability generating functions (PGFs) of the offspring distributions of the branching process, which determine respectively whether a major outbreak can occur and, if so, its probability.  The offspring distribution is different in the initial generation from that of all subsequent generations, since
the initial infective is chosen uniformly at random from the population (so its local and global degrees are independent), while subsequent primary infectives are infected through the network and their local and global degrees are dependent.
We focus first on the offspring means for a non-initial generation, since they determine whether or not the branching process is supercritical.

\subsubsection{Offspring mean matrix and threshold parameter $R_*$}
\label{rstar}
Let $\mathcal{B}_F$ denote the above multitype forward branching process and let $\tilde{\mathcal{B}}_F$ be the multitype branching process describing the descendants of a typical first-generation individual in $\mathcal{B}_F$.  Thus the type-dependent offspring law is the same for \emph{all} generations in
$\tilde{\mathcal{B}}_F$.
For $i=1,2,\cdots,n_Q$, let $\tilde{\bs{C}}_{i}=(\tilde{C}_{i1},\tilde{C}_{i2},
\cdots, \tilde{C}_{in_Q})$ be a vector random variable describing the numbers of offspring of different types of a typical type-$i$ individual in the branching process $\tilde{\mathcal{B}}_F$.  Thus,
$\tilde{C}_{ij}$ is the number of type-$j$ primary infectives generated by a typical single-household epidemic, whose primary infective is of type $i$.  Let $\tilde{M}=[\tilde{m}_{ij}]$ be the $n_Q\times n_Q$ matrix with elements
$\tilde{m}_{ij}={\rm E}[\tilde{C}_{ij}]$ and let $R_*$ be the dominant eigenvalue of $\tilde{M}$.  Then by standard multitype branching process theory
(see Mode (1971), Chapter 1, Theorem 7.1), the branching process $\tilde{\mathcal{B}}_F$ survives with strictly positive probability if and only if $R_*>1$.  Thus $R_*$ serves as a threshold parameter for our epidemic model. Note that this and subsequent results using the theory of multitype branching processes require assumptions regarding the irreducibility and/or positive regularity of the mean matrix $\tilde{M}$, which are met for all but highly pathological choices of $G$, $H$ and $n_Q$.

In order to compute $\tilde{M}$, and hence $R_*$, we need a further probability distribution.
For $d=1,2,\cdots$ and $h=1,2,\cdots,d$, let $\tilde{\pi}_h^{(d)}$ be the probability that a stub chosen uniformly at random from all stubs having total degree $d$ belongs to an individual who resides in a household of size $h$.  Note that this probability is the same for stubs with label $X=0$ and stubs with label $X=1$, and that
 \begin{equation*}
\tilde{\pi}_h^{(d)}=\frac{\pi_h h \tilde{p}_{d-h+1}}{\sum_{h'=1}^{d+1}\pi_{h'}h'\tilde{p}_{d-h'+1}}
=\frac{\tilde{\pi}_h \tilde{p}_{d-h+1}}{\sum_{h'=1}^{d+1}\tilde{\pi}_{h'}\tilde{p}_{d-h'+1}}.
\end{equation*}

To obtain $\tilde{m}_{ij}$, we condition first on the total degree of a typical type-$i$ primary infective and then on the size of its household yielding
\begin{equation}
\label{mtildeij}
\tilde{m}_{ij} = \sum_{d=1}^{\infty} p_{\tilde{D}|\tilde{Q}}(d|i) \sum_{h=1}^d \tilde{\pi}_h^{(d)}{\rm E}[\tilde{C}_{ij}^{(h,d)}],
\end{equation}
where $\tilde{\bs{C}}_{i}^{(h,d)}=(\tilde{C}_{i1}^{(h,d)},\tilde{C}_{i2}^{(h,d)},\cdots,\tilde{C}_{in_Q}^{(h,d)})$ is defined analogously to $\tilde{\bs{C}}_{i}$, except we condition on the type-$i$ individual residing in a household of size $h$ and having total degree $d$. (Note also that $p_{\tilde{D}|\tilde{Q}}(d|i)$ is is independent of the $X$-label of the individual concerned.)

Consider a typical size-$h$ single-household epidemic, with one initial infective, who is of type $i$ and has total degree $d$, and label the household members $0,1,\cdots,h-1$, where $0$ is the initial infective.  For $k=1,2,\cdots,h-1$, let $\chi_k=1$ if individual $l$ is infected by the single-household epidemic and let $\chi_k=0$ otherwise.  Then
\begin{equation}
\label{localdecomp}
\tilde{\bs{C}}_{i}^{(h,d)}=\tilde{\bs{C}}_{i}^{(h,d)}(0)+\sum_{k=1}^{h-1} \chi_k \tilde{\bs{C}}_{i}^{(h,d)}(k),
\end{equation}
where, for $k=0,1,\cdots,h-1$, $\tilde{\bs{C}}_{i}^{(h,d)}(k)=
(\tilde{C}_{i1}^{(h,d)}(k),\tilde{C}_{i2}^{(h,d)}(k),\cdots,\tilde{C}_{in_Q}^{(h,d)}(k))$, with $\tilde{C}_{ij}^{(h,d)}(k)$ being the number of type-$j$ primary infectives generated by individual $k$ in the single-household epidemic if it becomes infected.  (Throughout the paper, sums are zero if vacuous.)

Let $T^{(h)}=\sum_{k=1}^{h-1}\chi_k$ be the final size of the above single-household epidemic, not including the initial case, and let $\mu^{(h)}(\lambda)={\rm E}[T^{(h)}]$.  Then, see Ball (1986) equations (2.25) and (2.26),
\begin{equation*}
\mu^{(h)}(\lambda)=h-1-\sum_{k=0}^{h-1} \binom {h-1} {k} \alpha_k \phi_I(k\lambda)^{h-k}\quad (h=1,2,\cdots),
\end{equation*}
where $\phi_I(\theta)={\rm E}[\exp(-\theta I)]$ $(\theta \ge 0)$ is the moment generating function of $I$ and $\alpha_0,\alpha_1,\cdots$ are defined recursively by
\begin{equation*}
\sum_{l=0}^k \binom {k} {l} \alpha_l \phi_I(l\lambda)^{k-l}=k \quad (k=0,1,\cdots).
\end{equation*}

Note that $\chi_k$ and  $\tilde{\bs{C}}_{i}^{(h,d)}(k)$ are independent, because whether or not an individual is infected by the single-household epidemic is independent of its infectious period, so taking expectations of (\ref{localdecomp}) and noting that $\tilde{\bs{C}}_{i}^{(h,d)}(1),\tilde{\bs{C}}_{i}^{(h,d)}(2),\cdots, \tilde{\bs{C}}_{i}^{(h,d)}(h-1)$ are identically distributed yields
\begin{equation}
\label{mtildedecomp}
{\rm E}[\tilde{C}_{ij}^{(h,d)}]={\rm E}[\tilde{C}_{ij}^{(h,d)}(0)]+
\mu^{(h)}(\lambda){\rm E}[\tilde{C}_{ij}^{(h,d)}(1)].
\end{equation}

To determine ${\rm E}[\tilde{C}_{ij}^{(h,d)}(k)]$ $(k=0,1)$, for $i,j=1,2,\cdots,n_Q$ and $l=0,1$, let $p_{i,j}^{(l)}(r)$ be the probability that, when constructing the network, a given stub with $X$-label $l$ and total degree quantile $i$ is paired with a stub having total degree quantile $j$.  Then, $p_{i,j}^{(0)}=1/n_Q$ and
\begin{equation*}
p_{i,j}^{(1)}(r)=\begin{cases}
\delta_{i,j} & \mbox{if } r>0, \\
\delta_{i,n_Q+1-j} & \mbox{if } r<0,
\end{cases}
\end{equation*}
where $\delta_{i,j}=1$ if $i=j$ and $\delta_{i,j}=0$ if $i \ne j$.
Further, for $d=1,2,\cdots, j=1,2,\cdots,n_Q$ and $l=0,1$, let $\tilde{p}_{d,j}^{(l)}(r)$ be the probability that a stub chosen uniformly from all stubs having total degree $d$ and $X$-label $l$ is paired with a stub from quantile $j$.  Then $\tilde{p}_{d,j}^{(0)}(r)=1/n_Q$ and
\begin{equation*}
\tilde{p}_{d,j}^{(1)}(r)=\sum_{i=1}^{n_Q} p_{\tilde{Q}|\tilde{D}}(i|d)p_{i,j}^{(1)}(r).
\end{equation*}

Consider the individual labelled $0$, i.e.~the primary case, in the above single-household epidemic.  This individual has total degree $d$ and resides in a household of size $h$, so it has $d-h+1$ global neighbours, one of whom infected it.  Thus the individual has $d-h$ global edges along which it can spread the epidemic.  Each of the corresponding stubs  independently has $X$-label $1$ with probability $|r|$, so
\begin{equation}
\label{meanprim}
{\rm E}[\tilde{C}_{ij}^{(h,d)}(0)]=(d-h)p_I[(1-|r|)n_Q^{-1}+|r|p_{i,j}^{(1)}],
\end{equation}
where $p_I=1-\phi_I(\lambda)$ is the unconditional probability that a given infective infects a given susceptible neighbour.

Now consider the individual labelled $1$ in the single-household epidemic and suppose that it becomes infected.  The global degree of individual $1$ is distributed according to $G$.
Thus, for $g=1,2,\cdots$, with probability $p_g$, individual $1$ has $g$ global neighbours and hence total degree $g+h-1$.  Each of these
$g$ global neighbours is infected with probability $p_I$ and the $X$-labels of the corresponding outgoing stubs from individual $1$ are independent Bernoulli random variables with success probability $|r|$.  Summing over $g$ and taking expectations yields
\begin{equation}
\label{meansec}
{\rm E}[\tilde{C}_{ij}^{(h,d)}(1)]=\sum_{g=1}^{\infty}p_g g p_I[(1-|r|)n_Q^{-1}+|r|\tilde{p}_{g+h-1,j}^{(l)}(r)].
\end{equation}
Note that if $g=0$ then individual $1$ has no global neighbour to infect.  Note also that
${\rm E}[\tilde{C}_{ij}^{(h,d)}(1)]$ is independent of both $d$ and $i$, as indeed is the distribution of $\tilde{\bs{C}}_{i}^{(h,d)}(1)$.

Combining \eqref{mtildeij}, \eqref{mtildedecomp}, \eqref{meanprim} and \eqref{meansec} gives
\begin{align}
\tilde{m}_{ij} = p_I\sum_{d=1}^{\infty} p_{\tilde{D}|\tilde{Q}}(d|i) & \sum_{h=1}^d  \tilde{\pi}_h^{(d)} \biggl\{ (d-h)p_I \left[ (1-|r|)n_Q^{-1}+|r|p_{i,j}^{(1)} \right]  \nonumber \\
& +  \mu^{(h)}(\lambda) \biggl[ (1-|r|){\rm E}[G]n_Q^{-1} + |r|\sum_{g=1}^{\infty} p_g g \tilde{p}_{g+h-1,j}^{(l)}(r) \biggr] \biggr\} . \label{mean-matrix}
\end{align}
To summarise, equation (\ref{mean-matrix}) defines the elements of the mean matrix $\tilde M=[\tilde m_{ij}]$ of the branching process $\tilde{\mathcal{B}}_F$. The dominant eigenvalue of $\tilde M$, denoted by $R_*$, determines whether or not a major outbreak is possible, as described at the beginning of the section.

\subsubsection{Offspring PGFs and major outbreak probability}
\label{pglobal}

We now derive the offspring PGFs for the multitype branching processes $\mathcal{B}_F$ and $\tilde{\mathcal{B}}_F$, which enable their extinction probabilities (and hence the probability of a major outbreak) to be determined.
Observe that if the infectious periods are not constant, i.e.~there does not exist $\iota>0$ such that ${\rm P}(I=\iota)=1$, then the infectious periods of individuals infected by a single-household epidemic are not independent of the final size of that epidemic, which complicates, for example, using the decomposition~\eqref{localdecomp} to determine the offspring PGFs of $\tilde{\mathcal{B}}_F$.  As in Ball et al.~(2010), it is possible to use the theory of final state random variables developed in Ball and O'Neill (1999) to obtain expressions for these offspring PGFs in terms of Gontcharoff polynomials, though the details are rather involved and we do not present them here.  Instead, we consider the special case of a constant infection period, when the above-mentioned difficulties do not arise.  Thus in this subsection, but not elsewhere in Section~\ref{epid-prop}, we assume that $I\equiv\iota$ (i.e.~${\rm P}(I=\iota)=1$), so any given infective infects each of its neighbours (local or global) independently with probability $p_I=1-\exp(-\lambda \iota)$. The epidemic model is then an extension of the standard Reed-Frost epidemic (see, for example, Andersson and Britton (2000), Chapter 1) to our network model.  Note also that, in a physics setting, this Reed-Frost type model can be viewed as an extension, to incorporate degree correlation, of the bond percolation model of Gleeson (2009) for a class of clustered networks. Recall also that, as is well known for Reed-Frost type epidemics, the probability and the expected relative final size of a major outbreak are equal (cf.\ final paragraph of Section~\ref{zglobal}).

As noted previously, the forward branching process $\mathcal{B}_F$ has a different offspring distribution in the initial generation than in all subsequent generations.  We consider first a non-initial generation.  For $i=1,2,\cdots,n_Q$ and $\bs{s}=(s_1,s_2,\cdots,s_{n_Q})$ with $0 \le s_{i} \le 1$ $(i=1,2,\cdots,n_Q)$, let
\begin{equation*}
f_{\tilde{\bs{C}}_{i}}(\bs{s})={\rm E}\left[\prod_{j=1}^{n_Q}
s_{j}^{\tilde{C}_{ij}}\right]
\end{equation*}
be the joint PGF of $\tilde{\bs{C}}_{i}$.  (Throughout the paper, for a vector random variable, $\bs{Y}=(Y_1,Y_2,\cdots,Y_{n_Q})$ say, we use $f_{\bs{Y}}(\bs{s})$ to denote its joint PGF.)
Conditioning on the household size and total degree of a typical type-$i$ primary infective, as at (\ref{mtildeij}), yields
\begin{equation}
\label{PGFCitilde}
f_{\tilde{\bs{C}}_{i}}(\bs{s})=\sum_{d=1}^{\infty}p_{\tilde{D}|\tilde{Q}}(d|i)\sum_{h=1}^d f_{\tilde{\bs{C}}_{i}^{(h,d)}}(\bs{s}).
\end{equation}
The decomposition (\ref{localdecomp}) may be expressed as
\begin{equation}
\label{localdecomp1}
\tilde{\bs{C}}_{i}^{(h,d)}=\tilde{\bs{C}}_{i}^{(h,d)}(0)+\sum_{k=1}^{T^{(h)}} \tilde{\bs{C}}_{i}^{(h,d)}(k),
\end{equation}
where now $\tilde{\bs{C}}_{i}^{(h,d)}(1),\tilde{\bs{C}}_{i}^{(h,d)}(2),\cdots,
\tilde{\bs{C}}_{i}^{(h,d)}(T^{(h)})$ give the offspring vectors for the $T^{(h)}$
secondary cases in the single-household epidemic.  Further, since the infectious period is constant, conditional upon $T^{(h)}$, the random vectors $\tilde{\bs{C}}_{i}^{(h,d)}(1),\tilde{\bs{C}}_{i}^{(h,d)}(2),\cdots, \tilde{\bs{C}}_{i}^{(h,d)}(T^{(h)})$ are independent and identically distributed
copies of a random vector whose distribution is independent of $T^{(h)}$.  Hence, (\ref{localdecomp1}) implies that
\begin{equation}
\label{PGFdecomp1}
f_{\tilde{\bs{C}}_{i}^{(h,d)}}(\bs{s})=f_{\tilde{\bs{C}}_{i}^{(h,d)}(0)}(\bs{s})
f_{T^{(h)}}\left(f_{\tilde{\bs{C}}_{i}^{(h,d)}(1)}(\bs{s})\right),
\end{equation}
where $f_{T^{(h)}}(s)$ $(0 \le s \le 1)$ is the PGF of $T^{(h)}$, which, using Ball (1986), Theorem 2.6, is given by
\begin{equation}
\label{PGFTh1}
f_{T^{(h)}}(s)=s^{h-1}\sum_{k=0}^{h-1} \binom{h-1}{k}\alpha_k(s)(1-p_I)^{k(h-k)} \quad (h=1,2,\cdots),
\end{equation}
where $\alpha_0(s), \alpha_1(s),\cdots$ are defined recursively by
\begin{equation}
\label{PGFTh2}
\sum_{l=0}^k \binom{k}{l} (1-p_I)^{l(k-l)} \alpha_l(s)= s^{-k}
\quad (k=0,1,\cdots).
\end{equation}

To complete the derivation of $f_{\tilde{\bs{C}}_{i}}(\bs{s})$, we obtain expressions for $f_{\tilde{\bs{C}}_{i}^{(h,d)}(0)}(\bs{s})$ and
$f_{\tilde{\bs{C}}_{i}^{(h,d)}(1)}(\bs{s})$.
Consider a typical type-$i$ primary infective, $i^*$ say, and let $j^*$ be a susceptible global neighbour of $i^*$.  Let $\bs{\chi}_i=(\chi_{i1},\chi_{i2},\cdots,\chi_{in_Q})$, where
$\chi_{ik}=1$ if $i^*$ infects $j^*$ \emph{and} the edge between $i^*$ and $j^*$ was formed by connecting to a stub from $j^*$ belonging to quantile $k$, and $\chi_{ik}=0$
otherwise.  (Note that if $i^*$ does not infect $j^*$ then every element of  $\bs{\chi}_i $ is zero, and if $i^*$ does infect $j^*$ then precisely one element of $\bs{\chi}_i $ is one and all other elements of $\bs{\chi}_i $ are zero.)  For $i=1,2,\cdots,n_Q$ and $\bs{s}\in[0,1]^{n_Q}$, define the PGF of $\bs{\chi}_i$
\begin{equation}
g_i(\bs{s})={\rm E}\left[\prod_{j=1}^{n_Q}
s_{j}^{\chi_{ij}}\right]=1-p_I+p_I\sum_{j=1}^{n_Q}
\left[(1-|r|)\frac{s_j}{n_Q}+|r| p_{i,j}^{(1)}(r) s_{j}\right].
\label{gDef}
\end{equation}
Then using a similar argument to the derivation of (\ref{meanprim}) yields
\begin{equation}
\label{PGFhd0}
f_{\tilde{\bs{C}}_{i}^{(h,d)}(0)}(\bs{s})=\left(g_d(\bs{s})\right)^{d-h}.
\end{equation}
Now consider a typical individual, $\tilde{i}^*$ say, infected by a single-household epidemic and suppose that $\tilde{i}^*$ has total degree $d$.  Let $\tilde{j}^*$ be a susceptible global neighbour of $\tilde{i}^*$ and define $\tilde{\bs{\chi}}_d=(\tilde{\chi}_{d1},\tilde{\chi}_{d1},\cdots,\tilde{\chi}_{dn_Q})$ in the same way as $\bs{\chi}_i$ but with $i^*$ and $j^*$ replaced by $\tilde{i}^*$ and $\tilde{j}^*$, respectively.  Letting
\begin{equation}
\tilde{g}_d(\bs{s})={\rm E}\left[\prod_{j=1}^{n_Q}
s_{j}^{\tilde{\chi}_{ij}}\right]=1-p_I+p_I\sum_{j=1}^{n_Q}
\left[(1-|r|)\frac{s_j}{n_Q}+|r| \tilde{p}_{d,j}^{(1)}(r) s_{j}\right],
\label{gtildeDef}
\end{equation}
a similar argument to the derivation of \eqref{meansec} yields
\begin{equation}
\label{PGFhd1}
f_{\tilde{\bs{C}}_{i}^{(h,d)}(1)}(\bs{s})=\sum_{g=0}^{\infty}p_g \left(\tilde{g}_{g+h-1}(\bs{s})\right)^g.
\end{equation}
Combining~\eqref{PGFCitilde},~\eqref{PGFdecomp1},~\eqref{PGFhd0} and~\eqref{PGFhd1} gives the PGF of the offspring random variable $\tilde{\bs{C}}_{i}$ for a typical type-$i$ individual in $\tilde{\mathcal{B}}_F$.

Consider now the initial generation of the forward branching process $\mathcal{B}_F$.  Since the initial infective, $i^*$ say, in the epidemic is not infected through the network, the ancestor in $\mathcal{B}_F$ is not typed according to its total degree.
Let $\bs{C}=(C_1,C_2,\cdots,C_{n_Q})$ denote the offspring random variable for the ancestor in $\mathcal{B}_F$.  Then, conditioning on  $i^*$'s global degree and household size,
\begin{equation}
\label{PGFinit}
f_{\bs{C}}(\bs{s})=\sum_{g=0}^{\infty}\sum_{h=1}^{\infty}p_g \tilde{\pi}_h f_{\bs{C}^{(h,g+h-1)}}(\bs{s}), 
\end{equation}
where, for $h=1,2,\cdots$ and $d=h+1, h+2,\cdots$, ${\bs{C}}^{(h,d)}$ denotes the offspring random variable for the ancestor given that $i^*$ resides in a household of size $h$ and has total degree $d$.  Analogous to~\eqref{localdecomp1}, ${\bs{C}}^{(h,d)}$ admits the decomposition
\begin{equation}
\label{localdecompinit}
\bs{C}^{(h,d)}=\bs{C}^{(h,d)}(0)+\sum_{k=1}^{T^{(h)}}\bs{C}^{(h,d)}(k),
\end{equation}
whence, as at \eqref{PGFdecomp1},
\begin{equation}
\label{PGFinithd}
f_{\bs{C}^{(h,d)}}(\bs{s})=f_{{\bs{C}}^{(h,d)}(0)}(\bs{s})
f_{T^{(h)}}\left(f_{{\bs{C}}^{(h,d)}(1)}(\bs{s})\right).
\end{equation}
Now $\bs{C}^{(h,d)}(1)\overset{D}{=}\tilde{\bs{C}}^{(h,d)}(1)$, so $f_{{\bs{C}}^{(h,d)}(1)}(\bs{s})$ is given by the right hand side of~\eqref{PGFhd1}.
Note that if $i^*$ has household size $h$ and total degree $d$, then, since all of its $d-h+1$ global neighbours are susceptible, its offspring distribution is the same as that of a secondary infective having total degree $d$ in a single size-$h$ household epidemic.  Thus,
\begin{equation}
\label{PGFinithd0}
f_{{\bs{C}}^{(h,d)}(0)}(\bs{s})=\left(\tilde{g}_{d-h+1}(\bs{s})\right)^{d-h+1}.
\end{equation}
The offspring PGF $f_{\bs{C}}$ of the ancestor in $\mathcal{B}_F$ now follows using \eqref{PGFinit},~\eqref{PGFinithd},~\eqref{PGFhd1} and~\eqref{PGFinithd0}.

We now determine the probability of a major outbreak.  Suppose that $R_*>1$.  For $i=1,2,\cdots,n_Q$, let $\sigma_{i}$ be the probability that the branching process $\tilde{\mathcal{B}}_F$ goes extinct given that there is one ancestor whose type is $i$, and let
$\bs{\sigma}=(\sigma_{1},\sigma_{2},\cdots,\sigma_{n_Q})$.  Then, (see, for example, Mode (1971), Section 1.7.1), $\bs{\sigma}$ is the unique solution in $[0,1)^{n_Q}$ of the equations
\begin{equation}
\label{extinct}
f_{\tilde{\bs{C}}_{i}}(\bs{\sigma})=\sigma_{i} \quad(i=1,2,\cdots,n_Q).
\end{equation}
By conditioning on the number and type of offspring of the ancestor in $\mathcal{B}_F$, the probability that the branching process $\mathcal{B}_F$ survives (and hence the probability that a major outbreak occurs) is
\begin{equation}
\label{pmajor}
p_{\rm maj}=1-f_{\bs{C}}(\bs{\sigma}).
\end{equation}

\subsection{Final outcome of a major outbreak}
\label{zglobal}
We now consider the relative final size of a major outbreak.  The main tool that we use is the \emph{susceptibility set} (Ball (2000), Ball and Lyne (2001) and Ball and Neal (2002)), which we now define. Label the $n$ nodes (individuals) $1,2,\cdots,n$.  For $i=1,2,\cdots,n$, by sampling from the infectious period distribution and the Poisson processes describing when $i$ makes infectious contact with its neighbours, construct a (random) list of who $i$ would have infectious contact with if $i$ was to become infected.  Then construct a directed random graph, with nodes $1,2,\cdots,n$, in which for any pair of nodes $(i,j)$, with $i \ne j$, there is a directed edge from $i$ to $j$ if and only if $j$ is in $i$'s list.  For $i=1,2,\cdots,n$, the susceptibility set of node $i$ is set of all nodes $j$ from which there is a chain of directed edges to $i$ (including $i$ itself).

Observe that a node, $i$ say, is ultimately infected by the epidemic if and only if the initial infective belongs to $i$'s susceptibility set.  Suppose that the population size $n$ is large.  Then, as with the early stages of the epidemic, we can approximate the susceptibility set of a node, $i^*$ say, chosen uniformly at random from the population by a households-based multitype branching process.  We first consider $i^*$'s \emph{local} susceptibility set, i.e.~the set of nodes in $i^*$'s household from which there is a chain of within-household directed edges to $i^*$ (including $i^*$ itself). We next consider each member, $j^*$ say, of $i^*$'s local susceptibility set and determine which of $j^*$'s global neighbours have a directed edge joining them to $j^*$.  The set of all such global neighbours of $i^*$'s household form the first generation of the (backward) approximating branching process, with each such global neighbour, $k^*$ say, (generation-$1$ individual in the branching process) being typed by the quantile of the corresponding stub from $k^*$.  The process is then repeated in the obvious fashion to obtain the second generation of the backward branching process, and so on.
Denote this branching process by $\mathcal{B}_B$.  As with the forward branching process, the offspring law of $\mathcal{B}_B$ is different in the initial generation from that of all subsequent generations.
Let $\tilde{\mathcal{B}}_B$ be the multitype branching process describing the descendants of a typical first-generation individual in $\mathcal{B}_B$.

We conjecture that, subject to mild conditions on the household size and global degree distributions, the expected relative final size of a major outbreak converges to the survival probability of $\mathcal{B}_B$ as $n \to \infty$.  This is proved formally in Ball et al.~(2009) for the model with constant household size and no global degree correlation (i.e.~$r=0$); however, the proof in Ball et al.~(2009) is long and we do not attempt here to adapt it to the present model. Further, assuming the conjecture is true, the argument in Ball et al.~(2012) can be used to show that the relative final size of a major outbreak converges in probability to the survival probability of $\mathcal{B}_B$ as $n \to \infty$.  The proof in Ball et al.~(2012) is also quite long and we do not attempt to adapt it to the present model.  The numerical illustrations in Section~\ref{numexamples} (see Figure~\ref{fig:pmajCgce} and the surrounding commentary) support the above conjecture.

We determine now the offspring PGFs for $\mathcal{B}_B$ and $\tilde{\mathcal{B}}_B$.  We do not assume that the infectious periods are constant.  Let $\bs{B}=(B_1,B_2,\cdots,B_{n_Q})$ denote the offspring random variable for the ancestor in $\mathcal{B}_B$ and, for $i=1,2,\cdots,n_Q$, let
$\tilde{\bs{B}}_i=(\tilde{B}_{i1},\tilde{B}_{i2},\cdots,\tilde{B}_{in_Q})$ denote the offspring random variable for a typical type-$i$ individual in $\tilde{\mathcal{B}}_B$.

Consider $\tilde{\bs{B}}_i$ first.  Let $k^*$ be as above and assume it has type $i$.  Then arguing as at \eqref{PGFCitilde} yields
\begin{equation}
\label{PGFBitilde}
f_{\tilde{\bs{B}}_{i}}(\bs{s})=\sum_{d=1}^{\infty}p_{\tilde{D}|\tilde{Q}}(d|i)\sum_{h=1}^d f_{\tilde{\bs{B}}_{i}^{(h,d)}}(\bs{s}),
\end{equation}
where $\tilde{\bs{B}}_{i}^{(h,d)}$ denotes the corresponding offspring random variable when $k^*$ belongs to a household of size $h$ and has total degree $d$.
Let $M^{(h)}+1$ denote the size of a typical local susceptibility set in a household of size $h$. For $l=0,1$, let $\tilde{\bs{B}}_{i}^{(h,d)}(l)=(\tilde{B}_{i1}(l),\tilde{B}_{i2}(l),\cdots,\tilde{B}_{in_Q}(l))$, where  $\tilde{B}_{ij}(0)$ is the number of type-$j$ global neighbours of $k^*$ that would attempt to infect $k^*$ if they become infected and $\tilde{B}_{ij}(1)$ is defined similarly but for any other member of $k^*$'s local susceptibility set.
Then, noting that infectious global neighbours of an individual make infectious contact with that individual independently, each with probability $p_I$,
\begin{equation*}
f_{\tilde{\bs{B}}_{i}^{(h,d)}}(\bs{s})=f_{\tilde{\bs{B}}_{i}^{(h,d)}(0)}(\bs{s})
f_{M^{(h)}}\left(f_{\tilde{\bs{B}}_{i}^{(h,d)}(1)}(\bs{s})\right),
\end{equation*}
where, for $d=1,2,\cdots$ and $h=1,2,\cdots,d+1$,
\begin{equation*}
f_{\tilde{\bs{B}}_{i}^{(h,d)}(0)}(\bs{s})=\left(g_d(\bs{s})\right)^{d-h} \qquad
\mbox{and} \qquad
f_{\tilde{\bs{B}}_{i}^{(h,d)}(1)}(\bs{s})=\sum_{g=0}^{\infty}p_g \left(\tilde{g}_{g+h-1}(\bs{s})\right)^g
\end{equation*}
and $g_i(\bs{s})$ and $\tilde{g}_i(\bs{s})$ are defined by~\eqref{gDef} and~\eqref{gtildeDef}.

Turning to the PGF of $\bs{B}$, similar arguments to the above show that, in an obvious notation,
\begin{equation}
\label{PGFzinit}
f_{\bs{B}}(\bs{s})=\sum_{g=0}^{\infty}\sum_{h=1}^{\infty}p_g \tilde{\pi}_h
f_{{\bs{B}}^{(h,g+h-1)}(0)}(\bs{s})
f_{M^{(h)}}\left(f_{{\bs{B}}^{(h,g+h-1)}(1)}(\bs{s})\right),
\end{equation}
where, for $d=0,1,\cdots,$ and $h=1,2,\cdots,d+1$,
\begin{equation*}
f_{{\bs{B}}^{(h,d)}(0)}=\left(\tilde{g}_{d-h+1}(\bs{s})\right)^{d-h+1} \qquad \mbox{and} \qquad f_{{\bs{B}}^{(h,d)}(1)}=\sum_{g=0}^{\infty}p_g \left(\tilde{g}_{g+h-1}(\bs{s})\right)^g.
\end{equation*}

The probability mass function (and hence the PGF) of $M^{(h)}$ may be determined using the following result (see Ball and Neal~(2002), Lemma 3.1).
For $h=2,3,\cdots$,
\begin{equation*}
{\rm P}(M^{(h)}=k)=\binom{h-1}{k}\phi_I((k+1)\lambda)^{h-1-k}{\rm P}(M^{(k)}=k-1) \qquad (k=0,1,\cdots,h-1),
\end{equation*}
where
\begin{equation*}
\sum_{l=1}^k \binom{k-1}{l-1}\phi_I(l\lambda)^{k-l}{\rm P}(M^{(l)}=l-1)=1
\qquad (k=1,2,\cdots).
\end{equation*}

It is readily shown that ${\rm E}[M^{(h)}]={\rm E}[T^{(h)}]$ ($h=1,2,\cdots$), see Lemma 1 in the appendix of Ball et al.~(1997), using which it follows that $\tilde{\mathcal{B}}_B$ and $\tilde{\mathcal{B}}_F$ have the same offspring mean matrix.  Thus the branching process $\mathcal{B}_B$ survives if and only if $R_*>1$.
For $i=1,2,\cdots,n_Q$, let $\xi_{i}$ be the probability that the branching process $\tilde{\mathcal{B}}_F$ goes extinct given that there is one ancestor whose type is $i$, and let $\bs{\xi}=(\xi_{1},\xi_{2},\cdots,\xi_{n_Q})$. Then, if $R_*>1$, $\bs{\xi}$ is the unique solution in $[0,1)^{n_Q}$ of the equations
\begin{equation*}
f_{\tilde{\bs{B}}_{i}}(\bs{\xi})=\xi_{i} \quad(i=1,2,\cdots,n_Q)
\end{equation*}
and, for $n$ suitably large, the relative final size of a major outbreak, $z$ say, is given approximately by
\begin{equation}
\label{zfinal}
z=1-f_{\bs{B}}(\bs{\xi}).
\end{equation}

There does not appear to exist a similar recursive expression for the PGF $f_{M^{(h)}}(s)$ to that for $f_{T^{(h)}}(s)$ given by~\eqref{PGFTh1} and~\eqref{PGFTh2}, except when the infectious period is constant. In this case $M^{(h)}$ and $T^{(h)}$ have the same distribution, from which it easily follows (using the PGF formulae in the preceding sections) that $\pmaj=z$.

\section{Epidemics on rewired networks}
\label{rewired}
\subsection{Properties of epidemics}
\label{rewiredprop}
We now extend the results of the previous section to the model in which the edges in a fraction $p_{RW}$ of households are rewired.

Suppose first that $p_{RW}=1$, so all household edges are rewired.  The early stages of an epidemic in the rewired network may be approximated by a multitype branching process as in Section~\ref{forwardBP}, except now a local epidemic is the spread of disease along red edges alone, each having the same household size label.  Such local epidemics are realisations of the acquaintance model studied by Diekmann et al.~(1998) and a special case of a standard SIR epidemic on a configuration-model random network, see, for example, Newman (2002b).
Note that, if $n$ is large, the graph of red edges in the rewired network is locally tree-like.  For $h=2,3,\cdots$, let $\hat{\mathcal{E}}^{(h)}$ denote an SIR epidemic, with one initial infective, on a tree in which each node has degree $h-1$, with infectious period distributed according to $I$ and infection rate $\lambda$. Then for large $n$, a local epidemic in the rewired process may be approximated by $\hat{\mathcal{E}}^{(h)}$ and all the results of Sections~\ref{estab} and~\ref{zglobal} continue to hold provided the single-household final size and susceptibility set random variables $T^{(h)}$ and $M^{(h)}$ are replaced by their corresponding rewired counterparts defined on $\hat{\mathcal{E}}^{(h)}$, which we denote by $\hat{T}^{(h)}$ and $\hat{M}^{(h)}$.  As usual, the approximation of a local epidemic by $\hat{\mathcal{E}}^{(h)}$ can be made exact in the limit as $n \to \infty$ via a coupling argument.

Each individual in households of size $2$ have precisely one red stub, so when the corresponding red stubs are paired up such individuals are partitioned into households of size $2$ as before, whence $\hat{T}^{(2)}\overset{D}{=}T^{(2)}$ and
$\hat{M}^{(2)}\overset{D}{=}M^{(2)}$.  Fix $h\ge2$ and consider a typical local epidemic $\hat{\mathcal{E}}^{(h)}$.  The initial infective in $\hat{\mathcal{E}}^{(h)}$ has $h-1$ susceptible neighbours, while any
subsequent infective in the local epidemic has $h-2$ susceptible neighbours.
Any given infective infects any given susceptible neighbour with probability
$p_I=1-\phi_I(\lambda)$.  Thus in the (single-type) branching process, $\hat{\mathcal{B}}_F^{(h)}$ say, which gives the size of successive generations of infectives in $\hat{\mathcal{E}}^{(h)}$, the ancestor has offspring mean $(h-1)p_I$ and all subsequent individuals have offspring mean $(h-2)p_I$, whence
\begin{equation}
\label{rewiremuTh}
\hat{\mu}^{(h)}(\lambda)={\rm E}[\hat{T}^{(h)}]=\begin{cases}
(h-1)p_I[1-(h-2)p_I]^{-1} & \mbox{if } p_I<\frac{1}{h-2}, \\
\infty & \mbox{if } p_I\ge\frac{1}{h-2}.
\end{cases}
\end{equation}

Suppose now that $I\equiv\iota$, so any infective in $\hat{\mathcal{E}}^{(h)}$ infects each of its neighbours independently with probability $p_I$.  Then the offspring distribution of the ancestor in $\hat{\mathcal{B}}_F^{(h)}$ is ${\rm Bin}(h-1,p_I)$ and the offspring distribution of any subsequent individual is ${\rm Bin}(h-2,p_I)$, where  ${\rm Bin}(n,p)$ denotes a binomial distribution having $n$ trials and success probability $p$.  Standard branching process arguments then yield that, for $h=1,2,\cdots,$
\begin{equation}
\label{rewirefPGF}
f_{\hat{T}^{(h)}}(s)=\left(1-p_I+p_I\tilde{f}^{(h)}(s)\right)^{h-1} \quad (0 \le s \le 1),
\end{equation}
where $\tilde{f}^{(h)}(s)$ is the unique solution in $[0,1]$ of the equation
\begin{equation*}
\tilde{f}^{(h)}(s)=s\left(1-p_I+p_I\tilde{f}^{(h)}(s)\right)^{h-2},
\end{equation*}
cf.~equations (17) and (18) of Newman(2002b); note that $\tilde{f}^{(h)}(s)$ is the PGF of the total progeny of a typical non-ancestor in $\hat{\mathcal{B}}_F^{(h)}$.

Consider now the branching process, $\hat{\mathcal{B}}_B^{(h)}$ say, that describes on a generation basis a typical local susceptibility set associated with $\hat{\mathcal{E}}^{(h)}$ and return to the case of a general infectious period distribution.  It is easily seen that the offspring distributions of the ancestor and any subsequent individual in
$\hat{\mathcal{B}}_B^{(h)}$ are ${\rm Bin}(h-1,p_I)$ and ${\rm Bin}(h-2,p_I)$, respectively, where $p_I=\-\phi_I(\lambda)$, so $f_{\hat{M}^{(h)}}(s)$ is given by the right hand side of~\eqref{rewirefPGF}.

Finally we consider the case when the rewiring probability $p_{RW} \in (0,1)$.  Then, for example, the size $T^{(h)}(p_{RW})$ of a typical local epidemic corresponding to households having size $h$ is distributed according to $\hat{T}^{(h)}$, with probability $p_{RW}$, and to $T^{(h)}$, with probability $1-p_{RW}$.  Thus, ${\rm E}[T^{(h)}(p_{RW})]=(1-p_{RW})\mu^{(h)}(\lambda)+ p_{RW}\hat{\mu}^{(h)}(\lambda), f_{T^{(h)}(p_{RW})}(s)=(1-p_{RW})f_{T^{(h)}}(s)+p_{RW} f_{\hat{T}^{(h)}}(s)$ and $f_{M^{(h)}(p_{RW})}(s)=(1-p_{RW})f_{M^{(h)}}(s)+p_{RW} f_{\hat{M}^{(h)}}(s)$.  The threshold parameter $R_*$, probability of a major epidemic $p_{\rm maj}$ and relative final size of a major outbreak $z$ now follow by appropriate substitution into the results in Sections~\ref{estab} and~\ref{zglobal}.
\subsection{Effect of rewiring}
\label{rewiredeff}
We now examine the qualitative effect of rewiring on the probability and relative final size of a major outbreak.  For the model with $r=0$, constant infectious period and fixed household size (i.e.~${\rm P}(H=h)=1$ for some $h$), Gleeson et al.~(2010) use an analytic argument to show that the bond percolation threshold (the value of $p_I$ so that $R_*=1$) is larger for the model with full rewiring ($p_{RW}=1$) than for the model with no rewiring ($p_{RW}=0$).  Miller (2009) proves a similar result, again using an analytic argument, for an alternative model of random clustered networks, involving triangles, and also shows that the relative final size $z$ of a major outbreak is smaller for the fully rewired network than for the corresponding model without rewiring.  Here we employ a coupling argument, similar to that in, for example, Mollison (1977) and Ball (1983), to prove that for our model $R_*, p_{\rm maj}$ and $z$ are all increasing functions of the rewiring probability $p_{RW}$.  The coupling argument is both intuitive and powerful. It may be extended to the model of Gleeson et al.~(2010), without the restriction of a common household size, to the models of Miller (2009) and Newman (2009), and to the extension of the latter model proposed by Karrer and Newman (2010) that incorporates more general subgraphs than triangles.

For $h=1,2,\cdots$, let $\mathcal{E}^{(h)}$ denote the single size-$h$ household epidemic introduced in Section~\ref{rstar}, so $T^{(h)}$ is the final size of $\mathcal{E}^{(h)}$ not including the initial infective. For fixed $h\ge2$, a realisation of $\mathcal{E}^{(h)}$, viewed in generations of infectives, may be constructed from a realisation of $\hat{\mathcal{B}}_F^{(h)}$ as follows.  The ancestor of $\hat{\mathcal{B}}_F^{(h)}$ corresponds to the initial infective in $\mathcal{E}^{(h)}$.  The number of individuals, $Z_1$ say, in the first generation in $\hat{\mathcal{B}}_F^{(h)}$ (i.e.~the offspring of the ancestor) give the number of people directly infected by the initial infective in $\mathcal{E}^{(h)}$.  The individuals so infected are obtained by sampling $Z_1$ individuals uniformly at random without replacement from the $h-1$ individuals in the household excluding the initial infective.  The sampled individuals form the first generation of infectives in $\mathcal{E}^{(h)}$.  We now consider each first-generation individual in the branching process $\hat{\mathcal{B}}_F^{(h)}$ in turn.  The immediate offspring of such a first-generation individual give the number of people with which the corresponding infective in $\mathcal{E}^{(h)}$ makes infectious contact.  The people so contacted are obtained by sampling uniformly at random without replacement from the $h-1$ individuals in the household excluding the infective under consideration.  It is possible that a person so contacted has already been infected in $\mathcal{E}^{(h)}$, in which case the corresponding birth in $\hat{\mathcal{B}}_F^{(h)}$ and all of the descendants of that individual in $\hat{\mathcal{B}}_F^{(h)}$ are ignored in the construction of $\mathcal{E}^{(h)}$.  The construction of $\mathcal{E}^{(h)}$ continues in the obvious fashion and terminates when there is no infective remaining in the household.

Observe that by construction the size of the epidemic $\mathcal{E}^{(h)}$ is not larger than that the total progeny of the branching process $\hat{\mathcal{B}}_F^{(h)}$, so $\hat{T}^{(h)} \overset{st}{\ge} T^{(h)}$, where $\overset{st}{\ge}$ denotes stochastic ordering, whence $\hat{\mu}^{(h)}(\lambda) \ge \mu^{(h)}(\lambda)$ and $f_{\hat{T}^{(h)}}(s) \le f_{T^{(h)}}(s)$ $(0 \le s \le 1)$.  Moreover, provided $\lambda \mu_I>0$, these inequalities are strict for all $h \ge 3$ and all $s \in [0,1)$.  It follows that, if all other parameters are held fixed, the threshold parameter $R_*$ is an increasing function of the rewiring probability $p_{RW}$, as is the probability of a major outbreak $p_{\rm maj}$ (assuming that the infectious period is constant).  When the infectious period is not constant, the above coupling can be extended to include the global degrees of individuals in such a way that infectives in the household epidemic $\mathcal{E}^{(h)}$ have the same global degree and make the same global infectious contacts as the corresponding individuals in the branching process $\hat{\mathcal{B}}_F^{(h)}$, from which it follows that $p_{\rm maj}$ is increasing in $p_{RW}$.  Moreover, if ${\rm P}(H\ge3)>0$ and $\lambda \mu_I>0$ then both $R_*$ and $p_{\rm maj}$ are strictly increasing in $p_{RW}$.

Turning to the final outcome of a major outbreak, for fixed $h \ge 2$, we can construct a realisation of the local susceptibility set $\mathcal{S}^{(h)}$ say, of an individual, $i^*$ say, who resides in a household of size $h$, from a realisation of the branching process $\hat{\mathcal{B}}_B^{(h)}$ as follows.  The local susceptibility set of $i^*$ is constructed on a generation basis. The ancestor of
$\hat{\mathcal{B}}_B^{(h)}$ corresponds to the individual $i^*$.  The first generation of $\hat{\mathcal{B}}_B^{(h)}$ gives the number of individuals in $i^*$'s household who would make infectious contact with $i^*$ if they were to become infected; who these individuals (who form the first generation of $\mathcal{S}^{(h)}$) are is then determined by sampling without replacement as above.  We next consider in turn each member, $j^*$ say, of the first generation of $\mathcal{S}^{(h)}$ and determine which of those individuals not currently in $\mathcal{S}^{(h)}$ would join the susceptibility set of $i^*$ by virtue of making infectious contact with $j^*$.  Suppose that $j^*$ is the $k$th first-generation member of $\mathcal{S}^{(h)}$ to be considered in this fashion.  Then any individual not currently in $\mathcal{S}^{(h)}$ has failed to infect $k$ individuals, so the probability that it fails to infect $j^*$ is given by $p_F(k)=\phi_I((k+1)\lambda)/\phi_I(k\lambda)$.  Moreover, since such individuals are distinct, they each fail to infect $j^*$ independently with probability $p_F(k)$.  Let $p_F(0)= \phi_I(\lambda)$.  We now prove that, as one would expect on intuitive grounds, for any $\lambda>0$, $p_F(k)\ge p_F(0)$ ($k=1,2,\cdots$), with strict inequality unless $I \equiv \iota$ for some $\iota \ge 0$.

Define the function $\eta$ by $\eta(\theta)=\log \phi_I(\theta)$ $(\theta \ge 0)$.  Then $\eta$ is a convex function, since $\phi_I$ is a moment generating function, and
$\eta(0)=0$.  Thus, $\eta(\lambda)\le\frac{1}{k+1}\eta((k+1)\lambda)$ and $\eta(k\lambda)\le \frac{k}{k+1}\eta((k+1)\lambda)$, whence
\begin{equation}
\label{etacomp}
\eta(\lambda)+\eta(k\lambda)\le\eta((k+1)\lambda),
\end{equation}
which implies that $p_F(k)\ge p_F(0)$ ($k=1,2,\cdots$).  Moreover, if the infectious period random variable $I$ is not almost surely constant then $\eta$ is a strictly convex function, so, provided $\lambda>0$, the inequality in~\eqref{etacomp} is strict
and $p_F(k)> p_F(0)$ ($k=1,2,\cdots$).

In view of the above result, the individuals who join the susceptibility set $\mathcal{S}^{(h)}$ by virtue of making infectious contact with $j^*$ may be determined as follows.  Let $Z_{j^*}$ be the number of immediate offspring of the individual in $\hat{\mathcal{B}}_B^{(h)}$ that corresponds to $j^*$ and note that
$Z_{j^*}\sim{\rm Bin}(h-2,1-p_F(0))$.  Given $Z_{j^*}$, sample $\hat{Z}_{j^*}$ from
the binomial distribution ${\rm Bin}\left(Z_{j^*},\frac{1-p_F(k)}{1-p_F(0)}\right)$ and then
sample $Z_{j^*}$ individuals uniformly at random without replacement from the $h-1$ individuals in the household excluding $j^*$.  Any individual in this latter sample that is not currently in $\mathcal{S}^{(h)}$ is added to  $\mathcal{S}^{(h)}$.  This process is repeated for all $j^*$ belonging to the first generation of $\mathcal{S}^{(h)}$, thus yielding the second generation of $\mathcal{S}^{(h)}$, and so on.  Observe that, by construction, any individual in $\mathcal{S}^{(h)}$ has a corresponding individual in $\hat{\mathcal{B}}_B^{(h)}$, so $\hat{M}^{(h)} \overset{st}{\ge} M^{(h)}$, whence $f_{\hat{M}^{(h)}}(s) \le f_{M^{(h)}}(s)$ $(0 \le s \le 1)$, with strict inequality for $h \ge 3$ and $0 \le s <1$ provided $\lambda \mu_I>0$.  It follows that the relative final size $z$ of a major outbreak is  increasing in the rewiring probability $p_{RW}$, and strictly increasing if ${\rm P}(H\ge3)>0$ and $\lambda \mu_I>0$.

\section{Numerical examples}
\label{numexamples}
In this section we explore some properties of our network epidemic model numerically. We restrict our attention to the Reed-Frost type version of our model, i.e.\ we assume that $I\equiv \iota$ for some $\iota>0$, which implies that $\pmaj=z$, and rather than dealing explicitly with $I$ and the contact rate $\lambda$ we refer to the marginal infection probability $p_I=1-\exp(\lambda\iota)$.
Also, we use the notation $\Poi$ and $\Poi^+$ for global degree and household size distributions, as in Section~\ref{Sec-Tuning}.

\begin{figure}
\begin{center}
%\resizebox{10cm}{!}{\includegraphics{pMajCgcePoin1kN1k}}
\resizebox{\hfigwidth}{!}{\includegraphics{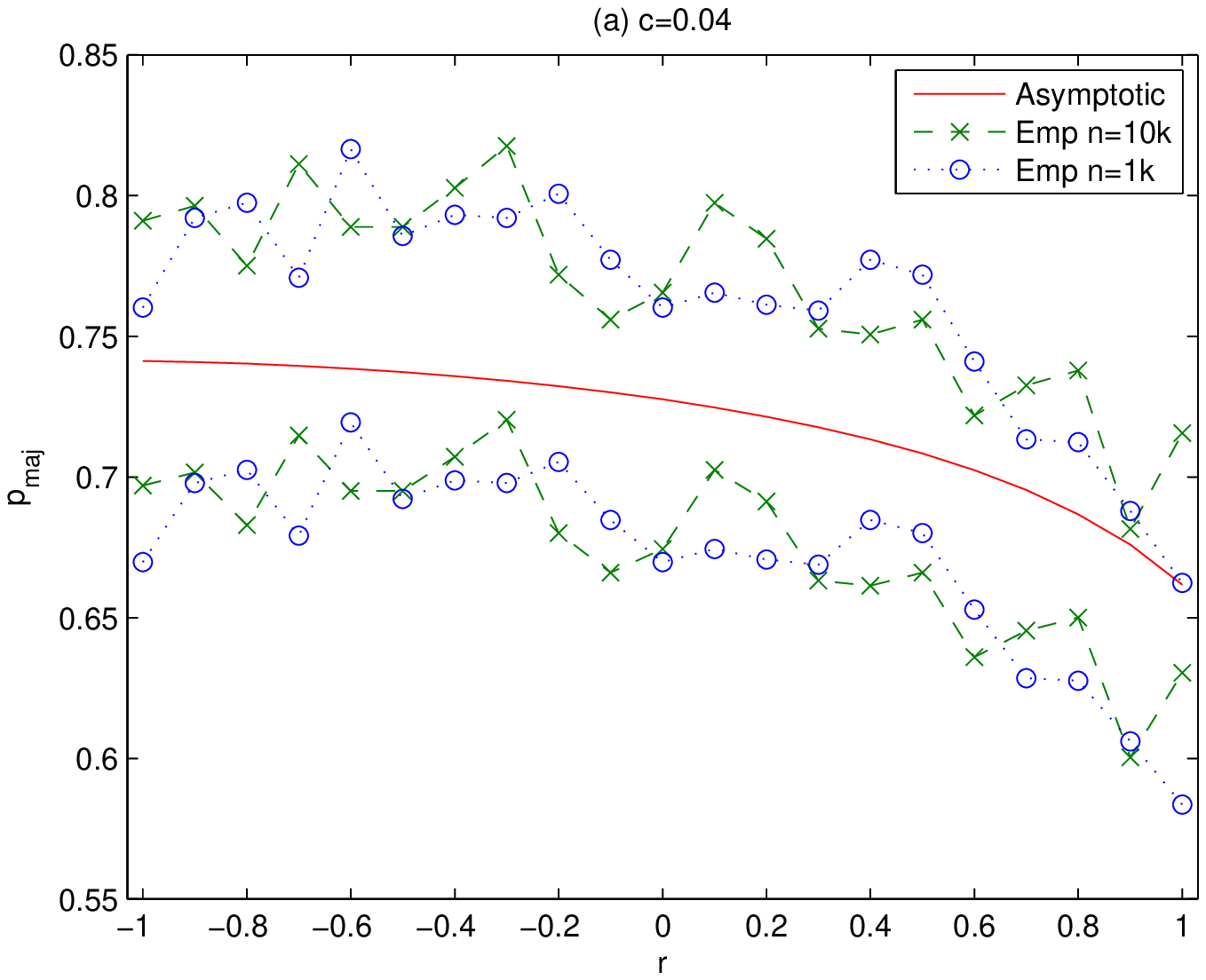}}
\resizebox{\hfigwidth}{!}{\includegraphics{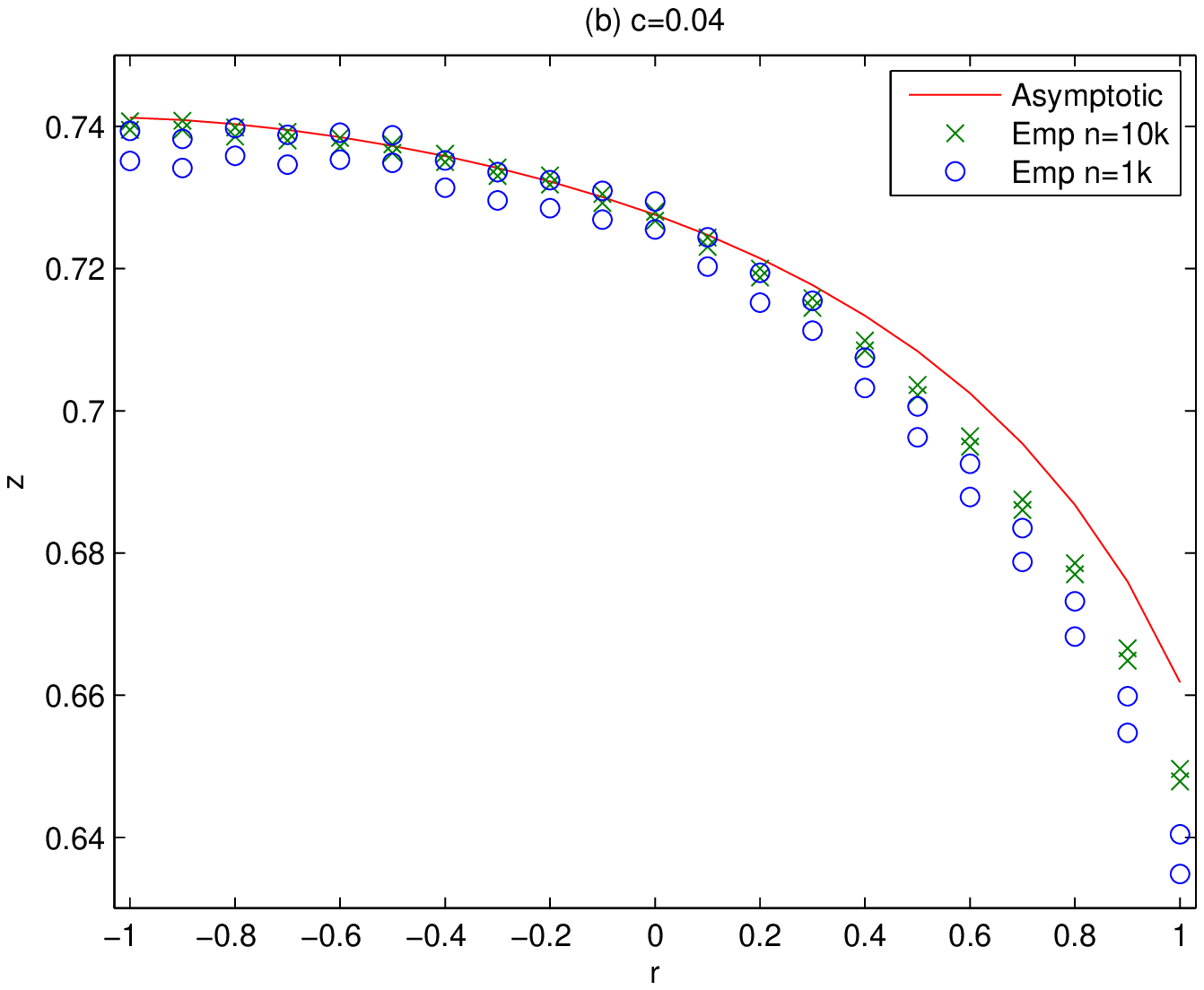}} \\
\resizebox{\hfigwidth}{!}{\includegraphics{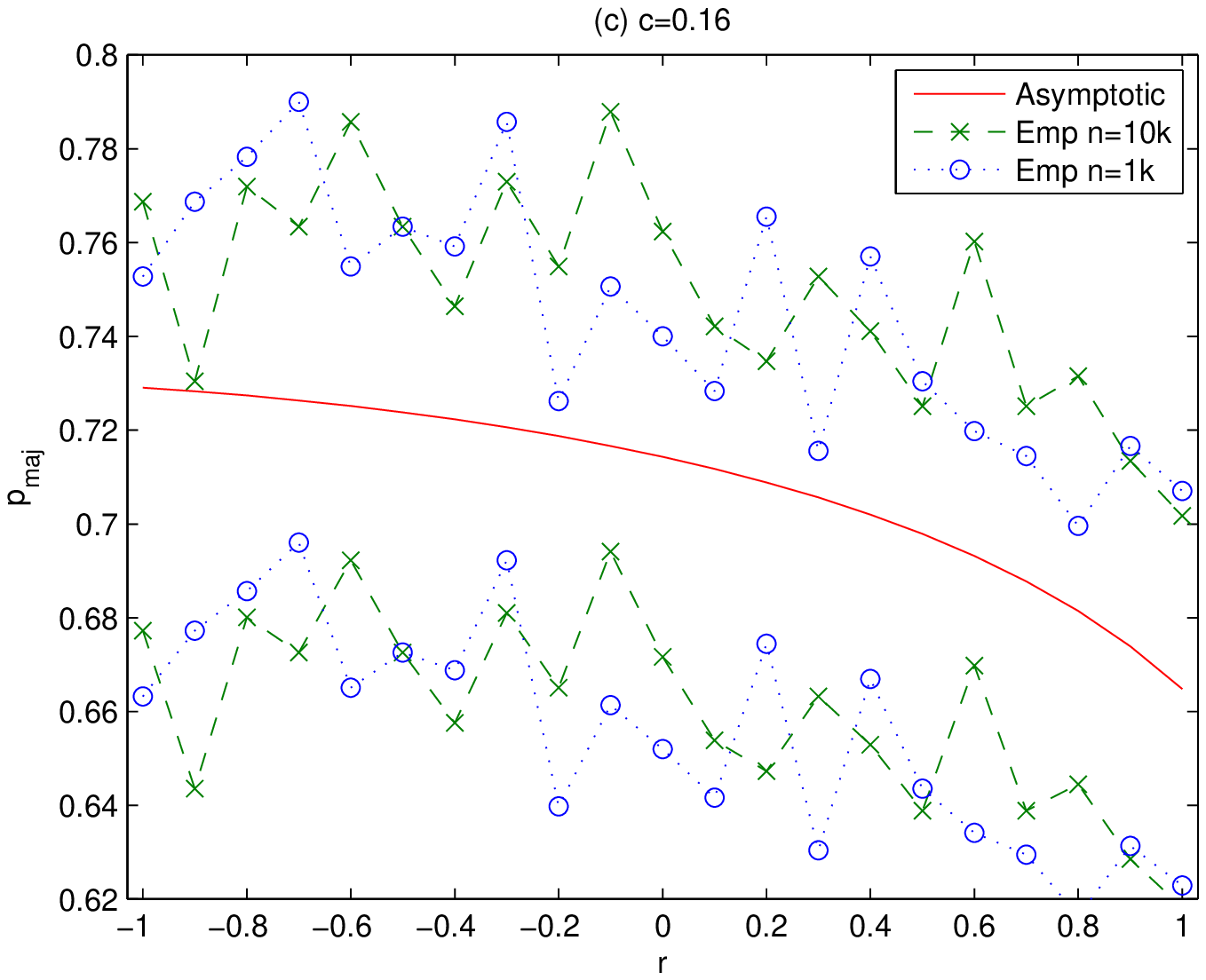}}
\resizebox{\hfigwidth}{!}{\includegraphics{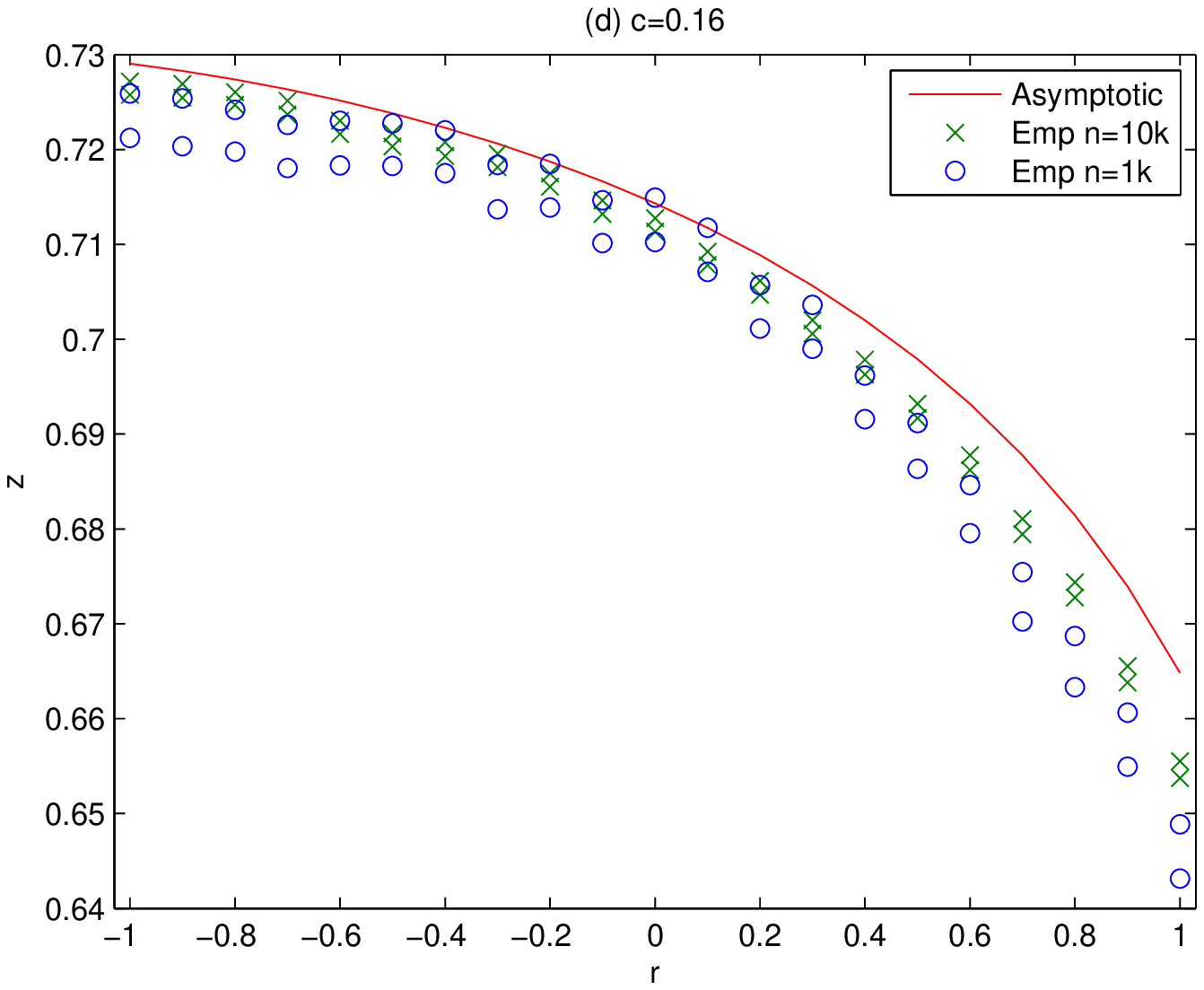}}
\end{center}
\caption{Plots comparing empirical estimates ($n<\infty$) and asymptotic values ($n\to\infty$) of $\pmaj$ and $z$, as a function of $r$, for our model with degree distributions $H\sim\Poi^+(2)$ and $G\sim\Poi(8)$ ($c=0.04$) and $H\sim\Poi^+(4)$ and $G\sim\Poi(6)$ ($c=0.16$). Other parameters are $n_Q=10$ and $p_I=0.2$. Empirical estimates are for network sizes $n=1,000$ and $n=10,000$, each estimate being based on 1,000 simulations. Note that the scales on the vertical axis on these plots is very variable.}
\label{fig:pmajCgce}
\end{figure}

First we briefly investigate the convergence of $\pmaj$ and $z$ for finite populations (derived empirically from simulations) to the asymptotic values (derived analytically) as the number of nodes/individuals $n$ becomes large. Figure~\ref{fig:pmajCgce} shows this behaviour in $\pmaj$ and $z$, for fixed $G$, $H$, $n_Q$, $p_I$ and varying $r\in[-1,1]$, comparing the asymptotic results to empirical estimates from networks of size $n=1,\!000$ and $10,\!000$ nodes/individuals. Each empirical estimate of a quantity of interest is based on $n_0=1,000$ simulations and is represented by an approximate 95.4\% confidence interval, calculated as a point estimate $\pm$ 2 standard errors (SE). (Also note that each simulation consists of generating a network then running an epidemic on it; we do not just run 1,000 epidemics on a single randomly generated network.) Each point estimate of $\pmaj$ is simply the proportion $\hat{p}$ of simulations that took off into a major outbreak (the cutoff between minor and major outbreaks being determined by inspecting histograms of epidemic final size), and ${\rm SE} = (\hat{p}(1-\hat{p})/n_0)^{1/2}$. The point estimate of $z$ is the mean fraction of the population ultimately infected by a major outbreak and here ${\rm SE} = \hat{\sigma} n_1^{-1/2}$, where $\hat{\sigma}^2$ is the sample variance of the  fraction of the population ultimately infected by a major outbreak and $n_1$ is the number of simulations that resulted in a major outbreak. As was explained in the closing sentences of Section 5 of Ball et al.\ (2009), our simulation methods yield much tighter confidence bands for $z$ than for $\pmaj$ since each simulation effectively gives a single realisation of the epidemic process but each simulation that results in a major outbreak gives $n-1$ (highly correlated) realisations of the susceptibility set process.%\marginal{Fig 2: Obviously would be better to have the higher clustering as 0.3--0.4 rather than 0.16, but there's not time and this gets the point across I think.}

We see that for networks with only 1000 nodes the asymptotic values of $\pmaj$ seem to be very good approximations to the empirically calculated major outbreak probabilities across all values of $r$. The expected relative final size also seems to be well approximated by the asymptotic values even for $n=1,000$; though there does appear to be some bias, which is more pronounced for more extreme values of $r$. One explanation for this is that when $r$ is close to $-1$ or $1$, there are more imperfections in the random graph (self-loops, household self-loops, etc.) and so the branching process approximation breaks down sooner. Nevertheless, the $z$ plots lend considerable credence to our conjecture in Section~\ref{zglobal} that the expected relative final size of a major outbreak converges to the survival probability of $\mathcal{B}_B$ as $n\to\infty$.

%{\bf Notes:} An additional difficulty here is that $n=1000$ is not sufficiently large for the cutoff (in histograms of simulated final sizes) between minor and major outbreaks to be clear. Taking $n$ up to 10,000 avoids this problem, probably 5000 would be sufficient though. Another alternative is of course to make the distinction clearer by using parameters that give a larger final size for major outbreaks; such simulations, however, take much longer to run because (i) more of them take off and (ii) those that do are larger. The simulations to get this Figure took about 12 hours. {\bf end notes}

\begin{figure}
\begin{center}
\resizebox{\hfigwidth}{!}{\includegraphics{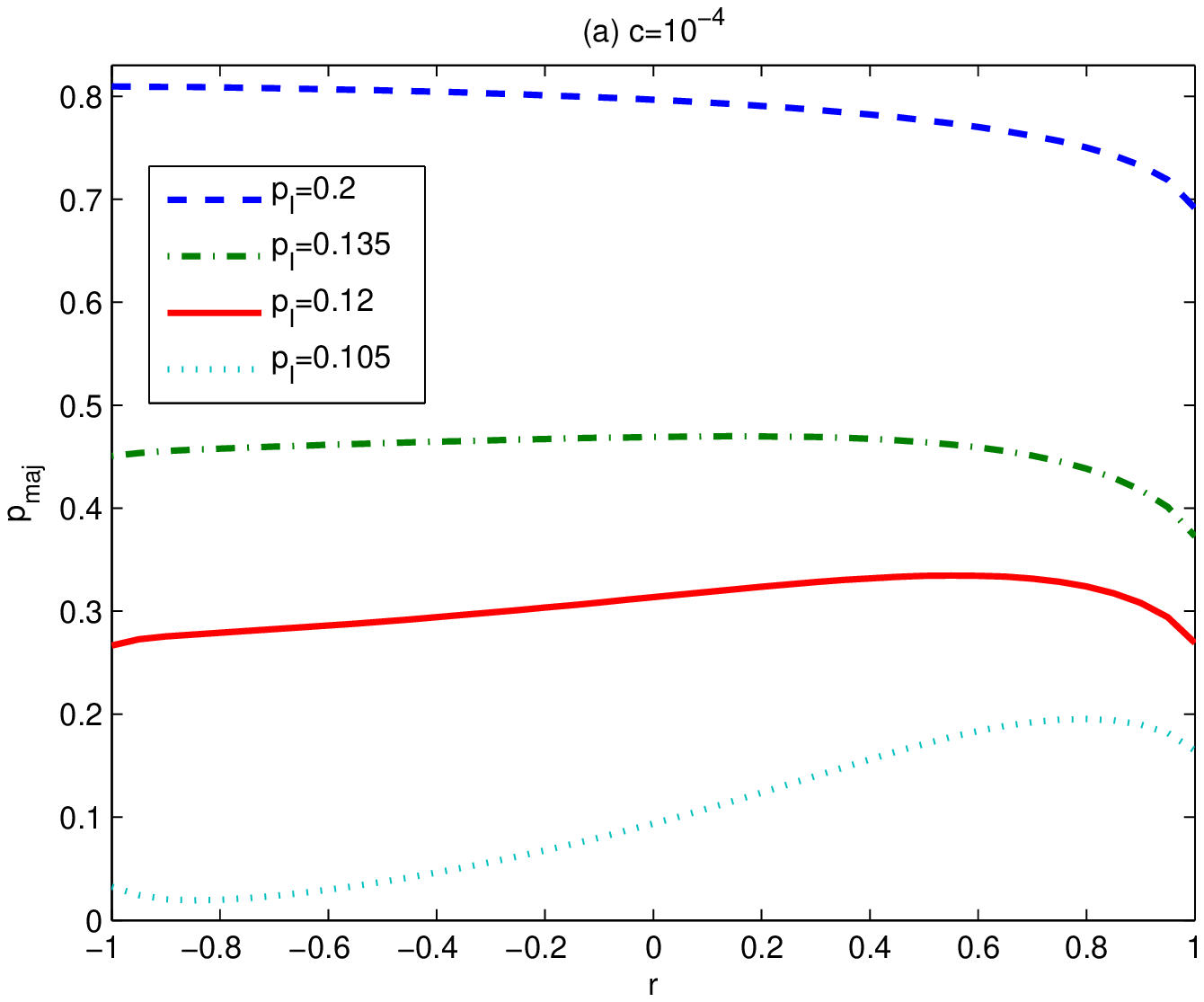}}
\resizebox{\hfigwidth}{!}{\includegraphics{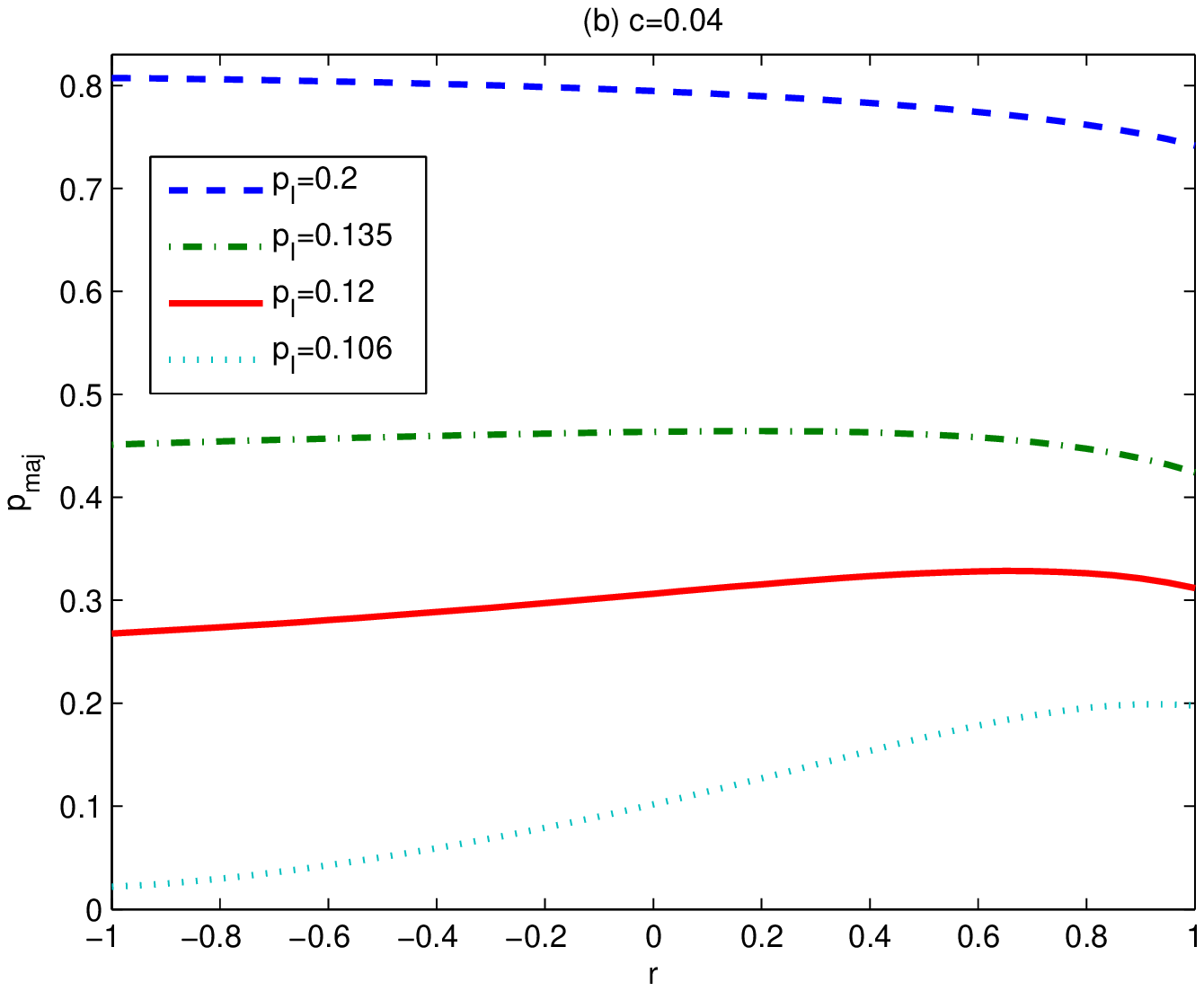}}\\
\resizebox{\hfigwidth}{!}{\includegraphics{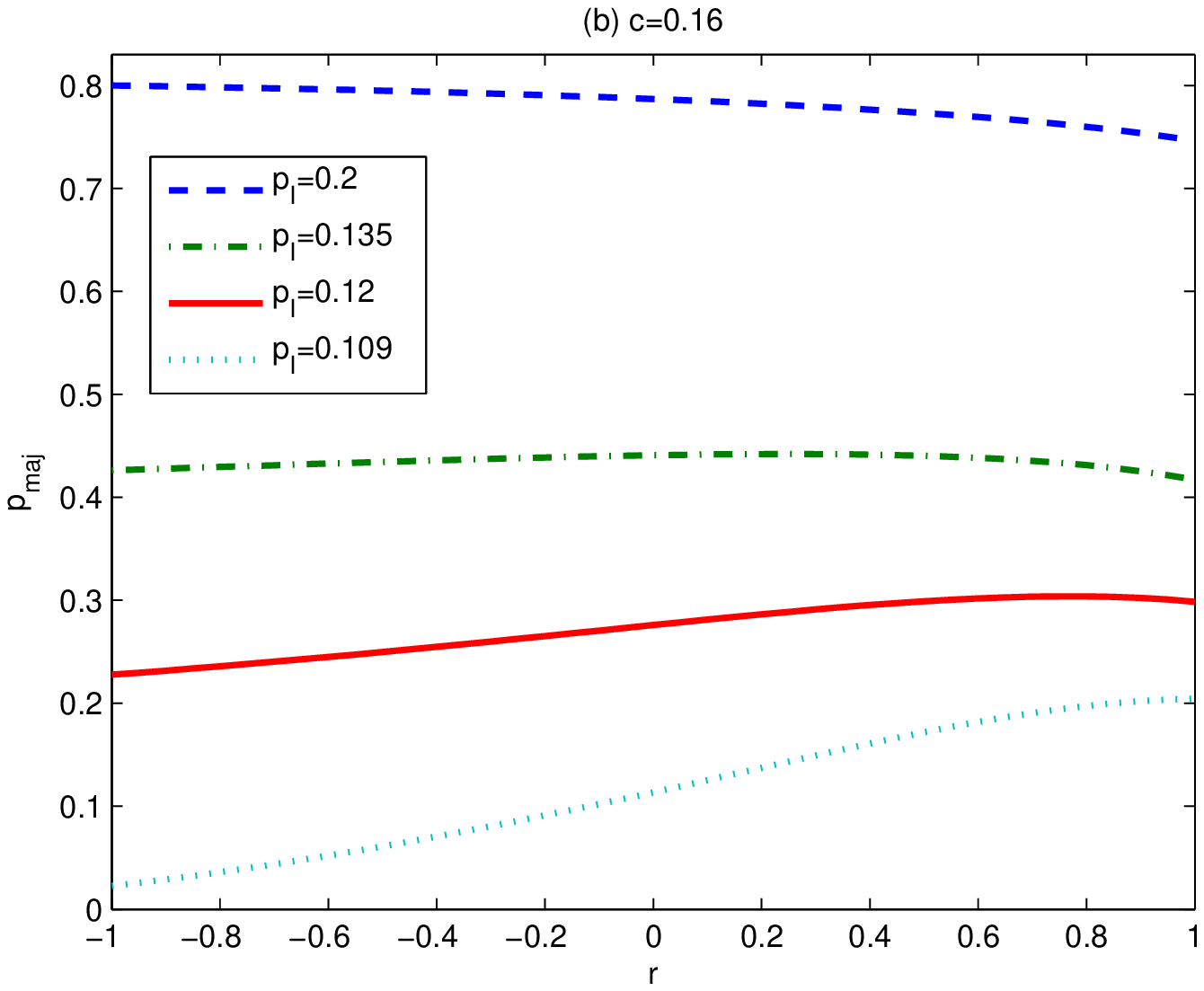}}
\resizebox{\hfigwidth}{!}{\includegraphics{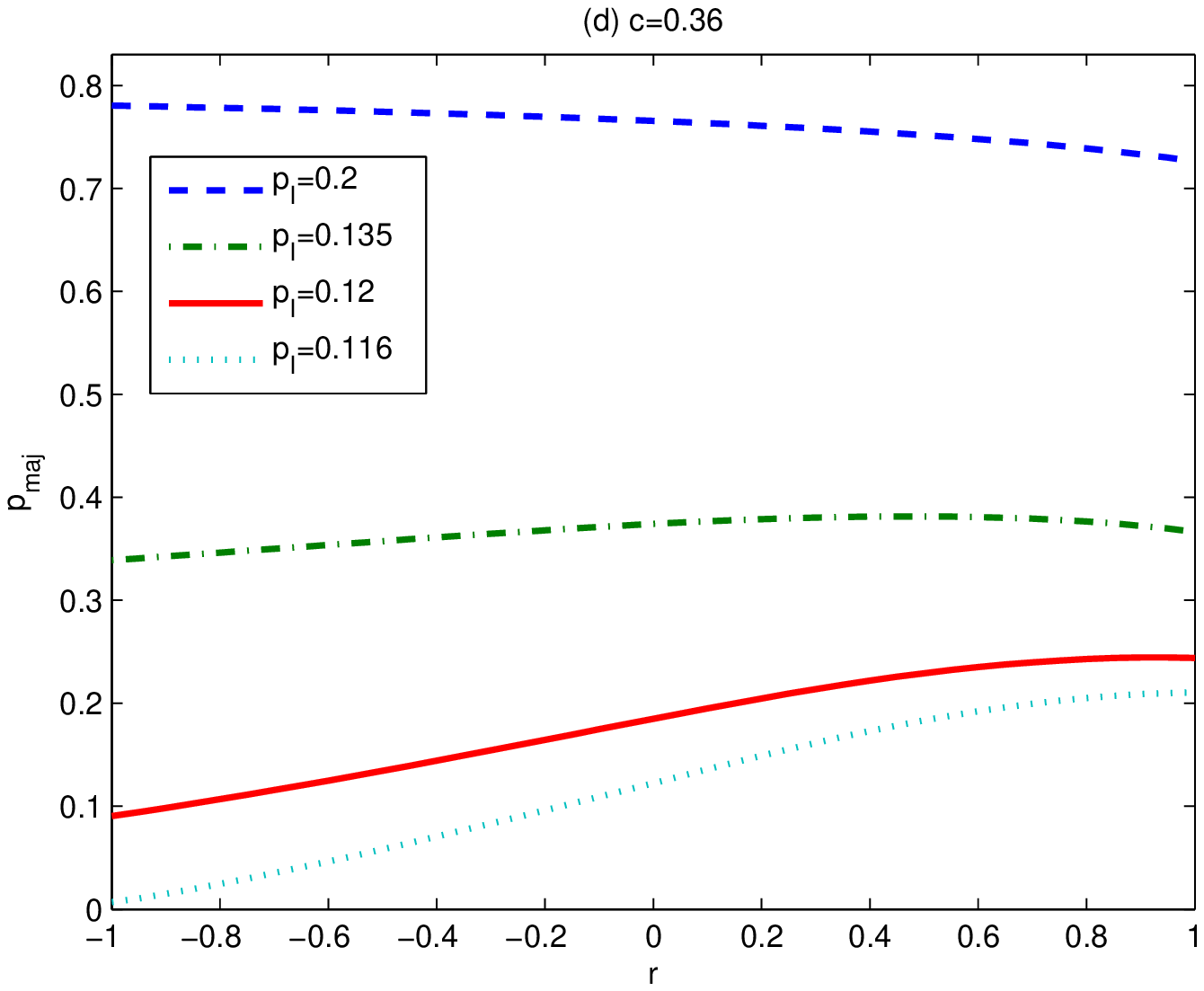}}
%\resizebox{10cm}{!}{\includegraphics{pMajVrPoiLarge.eps}}\\
%\resizebox{7cm}{!}{\includegraphics{pMajVrPoiG8H2Large.eps}}
%\resizebox{7cm}{!}{\includegraphics{pMajVrPoiG6H4Large.eps}}\\
%\resizebox{7cm}{!}{\includegraphics{pMajVrPoiG4H6Large.eps}}
%\resizebox{7cm}{!}{\includegraphics{pMajVrPoiG2H8Large.eps}}
\end{center}
\caption{Plot of $\pmaj$ versus $r$ for varying values of $p_I$. $G\sim\Poi(10-\mu)$ and $H\sim\Poi^+(\mu)$, with $\mu$ taking the values, in order, $0.1$, $2$, $4$, $6$; corresponding to clustering coefficients $10^{-4}, 0.04, 0.16, 0.36$. Note also that the $p_I$ values used are the same in each plot except for the smallest value, which is chosen so that the epidemic is just supercritical for all values of $r$.} 
\label{fig:pMajvRPoiLarge}
\end{figure}

Having seen that our asymptotic results give reasonable descriptions of the behaviour of our epidemic model on a moderately sized finite network, we turn our attention to investigating the effect of some of the parameters of our model on its (asymptotic) behaviour. We focus  initially on the qualitative behaviour of $\pmaj(=z)$ considered as a function of $r$ (and $p_I$). Figure~\ref{fig:pMajvRPoiLarge} illustrates this behaviour in the case where $G\sim\Poi(10-\mu)$, $H\sim\Poi^+(\mu)$, so $D\sim\Poi(10)$, and $n_Q=10$, for various values of $\mu\in[0,10)$ (and therefore $c=(\mu/10)^2$).
%\marginal{Fig 3: Choose 4 of these plots, ensure that the 4 lines on each one are in similar positions on the plot, and add headings `(a) $\mu=x$', etc.}

\begin{figure}
\begin{center}
\resizebox{\figwidth}{!}{\includegraphics{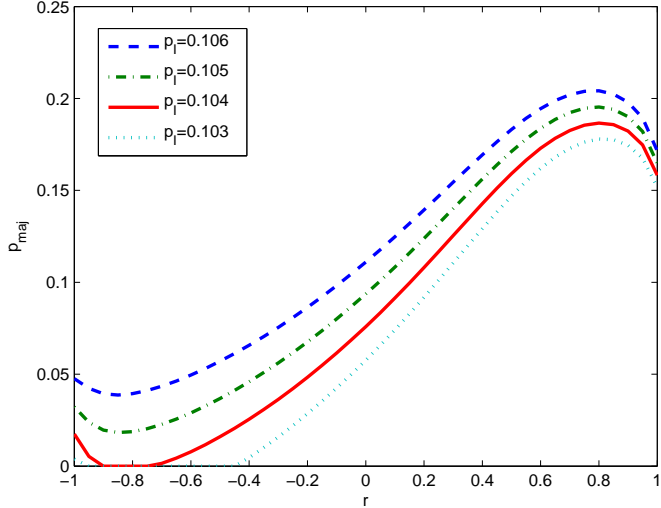}}
\end{center}
\caption{Plot of $\pmaj$ versus $r$ for near-critical values of $p_I$, when $G\sim\Poi(9.9)$, $H\sim\Poi^+(0.1)$ and $n_Q=10$. (Note that the $p_I=0.103$ line is positive near $r=-1$.)}
\label{fig:pMajvRPoiSmall}
\end{figure}

We see a variety of patterns in the dependance of $\pmaj$ on $r$ as $p_I$ and $c$ are varied. Broadly, when the process is well above criticality the dependance is not very strong, but when the process is only just supercritical changes in $r$ in particular (and thus in the degree correlation) can have a substantial impact on the epidemic model. The interesting (and somewhat unexpected) qualitative behaviour observed in the $p=0.105$ line in plot (a) is explored in further detail in Figure~\ref{fig:pMajvRPoiSmall}. Note, however, that the model parameters that give rise to this behaviour are $\mu_G=9.9$ and $\mu_H=0.1$, so there is essentially no clustering in the network; clearly further work is required to determine whether the model behaves in such a way with other, more realistic parameter values. Nevertheless, the wide range of values of $\pmaj(=z)$ for different values of $r$ (i.e.\ degree correlation) are observed near criticality in all of the plots in Figure~\ref{fig:pMajvRPoiLarge}; even though the non-monotonicity is only observed in plot~(a).

%{\bf Notes:} Obviously would be nicer to have this kind of behaviour when either (i) $H$ is less concentrated on 1 and/or (ii) $G$ follows some kind of heavy tailed distribution. I can get quite a bit of sensitivity to $r$ and $p_I$ using $H\sim\Poi(2)+1$ and $G\sim{\rm Pow}(5,3.2)$ ($\mu_G\approx 4.3$, $\sigma^2_G\approx 61.0$). Here $G\sim{\rm Pow}(k_*,a)$ means that $G$ has mass function
%\begin{equation*}
%p_k \propto \begin{cases}
%k_*^{-a}, & \mbox{for $k=1,2,\ldots,k_*$}, \\
%k^{-a}, & \mbox{for $k=k_*+1, k_*+2,\ldots$,}
%\end{cases}
%\end{equation*}
%which has a power law tail with index $a$. So far I can only get behaviour roughly the shape of the solid red line in Fig 3; without the decreasing bit near $r=-1$, but with the fairly substantial change from $r\approx-0.8$ to $r=1$; more numerical playing might reveal something more complex, but finding whether this is the case will take some time as the numerics are very slow when dealing with these heavy tailed distributions. {\bf end notes}

\begin{figure}
\begin{center}
\resizebox{\figwidth}{!}{\includegraphics{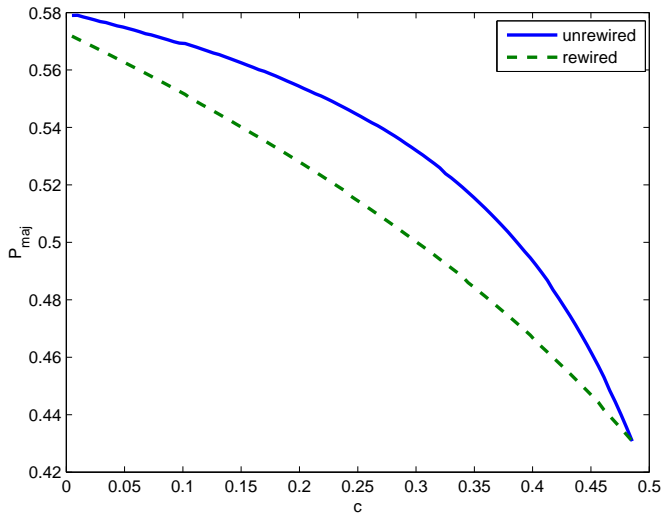}}
\end{center}
\caption{Plot of $\pmaj(=z)$ versus $c$ when $D\sim\Poi(10), \rho=0.2, n_Q=10$ and $p_I=0.15$.}
\label{fig:pMajvCPoi}
\end{figure}

Finally, Figure~\ref{fig:pMajvCPoi} illustrates the effect on $\pmaj(=z)$ of changing $c$, keeping $r$ and $p_I$ fixed, for the case when the total degree $D\sim\Poi(10)$ and $n_Q=10$.  The degree correlation $\rho$ is held fixed at $\rho=0.2$ and, for the unrewired model, the clustering coefficient $c$ is tuned to be any value in its feasible range (see Figure~\ref{fig:poissontune}) by varying $\mu$ and using~\eqref{Poicrho}.  The maximum value of $c$, consistent with $\rho=0.2$, is $c=0.4855$, which is attained when $r=-1$ and $\mu=6.9676$.  For the rewired model, the clustering coefficient is tuned by taking the unrewired model with $r=-1$ and $\mu=6.9676$ and letting the rewiring probability $p_{RW}$ vary in $[0,1]$.  
Figure~\ref{fig:pMajvCPoi} shows how $\pmaj(=z)$ varies with $c$ for both the unrewired and rewired models.  Note that, as one might expect, $\pmaj(=z)$ decreases with $c$ for both models; indeed this is proved formally for the rewired model in Section~\ref{rewiredeff}. Note also that $\pmaj(=z)$ is different for the two models, illustrating that these epidemic properties depend on more than just the local properties of the network encapsulated in $(D,c,r)$.  

\section{Discussion}
\label{sec-disc}

In this paper we define a network model which allows for quite arbitrary clustering $c$, degree correlation $\rho$ and degree distribution $D$, and asymptotic features of the model are derived. The main focus is on analysing an epidemic model on the network, and in particular what effect various network properties have on the epidemic in terms of its threshold parameter $R_*$, the probability $\pmaj$ of a major outbreak, and the relative size $z$ of a major outbreak. The main conclusion is that all three quantities $R_*$, $\pmaj$ and $z$ are decreasing with the clustering coefficient $c$ (when rewiring edges in the network thus keeping everything else fixed), whereas the dependence on the degree correlation $\rho$ is not as easily expressed: the quantities may be either increasing or decreasing depending on which part of the parameter space is being investigated. To our knowledge this is the first network model having such general features for which the properties of an epidemic are analysed in this level of detail.

A disadvantage with the model is that, in general, there is no simple and explicit relation between the model parameters $H$, $G$, $r$, and $n_Q$ and the more interesting network properties $c$, $\rho$ and $D$. Note however the relation for $D$ given in equation~\eqref{deg-dist}, and the facts that $\rho$ is increasing with $r$ and $c$ is increasing in $H$ (in the sense that $c(G,H_1,r) \ge c(G,H_2,r)$ if $H_1 \overset{st}{\ge} H_2$), keeping other parameters fixed. A model having simpler relationships to the local network properties could be more easily interpreted and would hence be of interest.  The use of appropriate pairing of stubs to control degree correlation, as done in this paper, could be applied to other models of clustered networks, such as those in Newman (2009), Miller (2009) and Karrer and Newman (2010).

It is important to observe that, as illustrated in Figure~\ref{fig:pMajvCPoi}, there may be distinct network models having the same local network features $D$, $\rho$ and $c$ but still giving different properties of an epidemic, the latter being a global property.  In applications it is hence important to fit not only local properties of a network model to empirical network data, but also to study the definitions of the model and try to understand if the model mechanism seems to agree realistically with how the empirical network may have been constructed.

\ifthenelse{\plain=1}
{{\bf Acknowledgments}}
{\begin{acknowledgements}}
This research was supported by the UK Engineering and Physical Sciences Research Council, under research grant number EP/E038670/1 (FB and DS), and by the Swedish Research Council (TB).
\ifthenelse{\plain=1}
{}
{\end{acknowledgements}}

\appendix
\section*{Appendix: Derivation of degree correlation $\rho$}
In the appendix we derive the formula for the degree correlation $\rho$ for our model given in equation~\eqref{degcorr}.  Let $E$ denote an edge chosen uniformly at random from all edges in the network, and let $X_L$ and $X_R$ denote the total degrees of the nodes adjacent to $E$. Then $\rho={\rm corr}(X_L,X_G)$, i.e.~the correlation between
$X_L$ and $X_R$.  Let $I_G=1$ if $E$ is a global edge and $I_G=0$ if $E$ is a household edge, so ${\rm P}(I_G=1)=p_G=1-{\rm P}(I_G=0)$.  We determine first the probability $p_G$ that $E$ is a global edge.

Let $N_G$ and $N_H$ denote respectively the number of global and household edges in the network.  Then $\mu_{N_G}=\frac{n}{2}\mu_G$, since each stub contributes to half an edge, and $\mu_{N_H}=\frac{n}{2}\mu_{\tilde{H}-1}$, since the household size of an individual chosen unifomly at random from the population is distributed according to $\tilde{H}$ and if such an individual resides in a household of size $h$ it has $h-1$ household neighbours.  Letting $n\to \infty$ and using the strong law of large numbers shows that $p_G$ is given by~\eqref{frac-glob}.

Note that
\begin{equation}
\label{covXLXR}
{\rm cov}(X_L,X_R)={\rm E}[{\rm cov}(X_L,X_R|I_G)]+{\rm cov}({\rm E}[X_L|I_G],{\rm E}[X_R|I_G]).
\end{equation}
We calculate the two quantities on the right hand side of~\eqref{covXLXR} in turn.

Suppose that $I_G=0$, so $E$ is a household edge.  Then $X_L=H_E-1+G_L$ and $X_R=H_E-1+G_R$, where $H_E$ is the size of the household that contains the edge $E$, and $G_L$ and $G_R$ are the global degrees of the nodes adjacent to $E$.  Observe that $H_E$ is distributed as $\hat{H}$ and, since $I_G=0$, $G_L$ and $G_R$ are independent copies of $G$.
Thus,
\begin{equation}
\label{covXLXRIG0}
{\rm cov}(X_L,X_R|I_G=0)=\sigma_{\hat{H}}^2.
\end{equation}

Suppose that $I_G=1$, so $E$ is a global edge.  Let $Q_L$ and $Q_R$ be the total degree quantiles of the two stubs used to form the edge $E$.  Then, for $i,j=1,2,\cdots,n_Q$,
\begin{equation}
\label{PQLQR}
{\rm P}(Q_L=i, Q_R=j)=\begin{cases}
\frac{1-r}{n_Q^2}+\delta_{i,j}\frac{r}{k} & \mbox{if } r \ge 0, \\
\frac{1-|r|}{n_Q^2}+\delta_{i,n_Q+1-j}\frac{|r|}{k} & \mbox{if } r<0.
\end{cases}
\end{equation}
Now,
\begin{align}
\label{covXLXRIG1}
{\rm cov}(X_L,X_R|I_G=1)&={\rm E}[{\rm cov}(X_L,X_R|I_G=1,Q_L,Q_R)]\nonumber\\&\qquad +{\rm cov}({\rm E}[X_L|I_G=1,Q_L],{\rm E}[X_R|I_G=1,Q_R]).
\end{align}
Given $(Q_L,Q_R)$, the total degrees $X_L$ and $X_R$ are independent, so
\begin{equation}
\label{covXLXRIG1Q}
{\rm cov}(X_L,X_R|I_G=1,Q_L,Q_R)=0.
\end{equation}
Further, for $i=1,2,\cdots,n_Q$, ${\rm E}[X_L|I_G=1,Q_L=i]={\rm E}[X_L|I_G=1,Q_R=i]=\mu_{\tilde{D}}^{(i)}$ (see equation~\eqref{mudtildei}).
Using the distribution~\eqref{PQLQR} and noting that $\mu_{\tilde{D}}=n_Q^{-1}\sum_{i=1}^{n_Q}\mu_{\tilde{D}}^{(i)}$ yields
\begin{equation}
\label{covEXLEXRIGQ}
{\rm cov}({\rm E}[X_L|I_G=1,Q_L],{\rm E}[X_R|I_G=1,Q_R])=g_{\tilde{D},n_Q}(r),
\end{equation}
where $g_{\tilde{D},n_Q}(r)$ is defined at~\eqref{gtilded}.

Note that ${\rm P}(I_G=1)=p_G=1-{\rm P}(I_G=0)$.  Then, equations~\eqref{covXLXRIG0},~\eqref{covXLXRIG1},~\eqref{covXLXRIG1Q} and~\eqref{covEXLEXRIGQ} yield
\begin{equation}
\label{EcovXLXRIG}
{\rm E}[{\rm cov}(X_L,X_R|I_G)]=(1-p_G)\sigma_{\hat{H}}^2+p_G g_{\tilde{D},n_Q}(r).
\end{equation}

We turn now to the second quantity on the right hand side of~\eqref{covXLXR}.  Note that ${\rm E}[X_L|I_G]={\rm E}[X_R|I_G]$, so ${\rm cov}({\rm E}[X_L|I_G],{\rm E}[X_R|I_G])={\rm var}({\rm E}[X_L|I_G])$.  Suppose that $I_G=0$.  Then, in the above notation, $X_L=H_E-1+G_L$, where $G_L\overset{D}{=}G$.  Thus,
\begin{equation}
\label{EXLIG0}
{\rm E}[X_L|I_G=0]=\mu_{\hat{H}-1}+\mu_G.
\end{equation}
Suppose that $I_G=1$.  Then $X_L\overset{D}{=}\tilde{D}$ and recall that $\tilde{D}\overset{D}{=}\tilde{H}-1+\tilde{G}$.  Thus,
\begin{equation}
\label{EXLIG1}
{\rm E}[X_L|I_G=1]=\mu_{\tilde{H}-1}+\mu_{\tilde{G}}.
\end{equation}
Recalling that ${\rm P}(I_G=1)=p_G=1-{\rm P}(I_G=0)$ and that $\mu_{\tilde{G}}={\rm E}[G^2]/\mu_G$, equations~\eqref{EXLIG0} and~\eqref{EXLIG1} yield
\begin{equation}
\label{covEXLEXRIG}
{\rm cov}({\rm E}[X_L|I_G],{\rm E}[X_R|I_G])=
p_G(1-p_G)\left(\mu_{\hat{H}}-\mu_{\tilde{H}}-\frac{\sigma_G^2}{\mu_G}\right)^2.
\end{equation}
Combining equations~\eqref{covXLXR},~\eqref{EcovXLXRIG} and~\eqref{covEXLEXRIG} gives
\begin{equation}
\label{covXLXR1}
{\rm cov}(X_L,X_R)=(1-p_G)\sigma_{\hat{H}}^2+p_G g_{\tilde{D},n_Q}(r)+p_G(1-p_G)\left(\mu_{\hat{H}}-\mu_{\tilde{H}}-\frac{\sigma_G^2}{\mu_G}\right)^2.
\end{equation}

We now derive ${\rm var}(X_L)$.  First note that
\begin{equation}
\label{varXL}
{\rm var}(X_L)={\rm E}[{\rm var}(X_L|I_G)]+{\rm var}({\rm E}[X_L|I_G]).
\end{equation}
As above, if $I_G=0$ then $X_L=H_E-1+G_L$, where $H_E\overset{D}{=}\hat{H}$ and
$G_L\overset{D}{=}G$ are independent, so ${\rm var}(X_L|I_G=0)=\sigma_{\hat{H}}^2+\sigma_G^2$; and if $I_G=1$ then $X_L\overset{D}{=}\tilde{H}-1+\tilde{G}$, where $\tilde{H}$ and $\tilde{G}$ are independent, so
${\rm var}(X_L|I_G=1)=\sigma_{\tilde{H}}^2+\sigma_{\tilde{G}}^2$.  Hence,
\begin{equation*}
{\rm E}[{\rm var}(X_L|I_G)]=(1-p_G)\left(\sigma_{\hat{H}}^2+\sigma_G^2 \right)+p_G\left(\sigma_{\tilde{H}}^2+\sigma_{\tilde{G}}^2\right),
\end{equation*}
which on substituting into~\eqref{varXL}, recalling that ${\rm var}({\rm E}[X_L|I_G])={\rm cov}({\rm E}[X_L|I_G],{\rm E}[X_R|I_G])$ and using~\eqref{covEXLEXRIG}
yields
\begin{equation}
\label{varXL1}
{\rm var}(X_L)=(1-p_G)\left(\sigma_{\hat{H}}^2+\sigma_G^2 \right)+p_G\left(\sigma_{\tilde{H}}^2+\sigma_{\tilde{G}}^2\right)
+p_G(1-p_G)\left(\mu_{\hat{H}}-\mu_{\tilde{H}}-\frac{\sigma_G^2}{\mu_G}\right)^2.
\end{equation}

The expression~\eqref{degcorr} for the degree correlation $\rho$, given in Section~\ref{sec-corr}, follows from
equations~\eqref{covXLXR1} and~\eqref{varXL1}, since ${\rm var}(X_L)={\rm var}(X_R)$.

\section*{References}

\ifthenelse{\plain=1}
{}
{\indent}

%{\sc Anderson, R.M.\ and May, R.M.} (1991),
%{\it Infectious Diseases of Humans; dynamics and control,}
%Oxford University Press, Oxford.

{\sc Andersson, H.} (1999),
Epidemic models and social networks,
{\it The Mathematical Scientist} {\bf 24(2)} 128--147.

{\sc Andersson, H. and Britton, T.} (2000), Stochastic epidemic
models and their statistical analysis, {\it Springer Lecture Notes
in Statistics} {\bf 151}, New York: Springer Verlag.

{\sc Badham, J.\ and Stocker, R.} (2010),
The impact of network clustering and assortativity on epidemic behaviour,
{\it Theor.\ Pop.\ Biol.} {\bf 77} 71--75.

{\sc Ball, F.G.} (1983),
The threshold behaviour of epidemic models,
{\it J.\ Appl.\ Prob.} {\bf 20} 227--241.

{\sc Ball, F.G.} (1986),
A unified approach to the distribution of total size and total area under the trajectory of the infectives in epidemic models,
{\it Adv.\ Appl.\ Prob.} {\bf 18} 289--310.

{\sc Ball, F.G.} (2000),
Susceptibility sets and the final outcome of stochastic {SIR} epidemic models.
Research Report 00-09. Division of Statistics, School of Mathematical Sciences, University of Nottingham.

{\sc Ball, F.G.\ and Lyne, O.D.} (2001),
Stochastic multitype SIR epidemics among a population partitioned into households,
{\it Adv.\ Appl.\ Prob.} {\bf 33} 99--123.

{\sc Ball, F.G.; Mollison, D.\ and Scalia-Tomba, G.} (1997),
Epidemics with two levels of mixing,
{\it Ann.\ Appl.\ Prob.} {\bf 7} 46-89.

{\sc Ball, F.G.\ and Neal, P.} (2002),
A general model for stochastic SIR epidemics with two levels of mixing,
{\it Math.\ Biosci.} {\bf 180} 73--102.

{\sc Ball, F.G.\ and O'Neill, P.D.} (1999),
The distribution of general final state random variables for stochastic epidemic models,
{\it J.\ Appl.\ Prob.} {\bf 36} 473--491.

{\sc Ball, F.G.\ and Sirl, D.J.} (2012),
An SIR epidemic model on a population with random network and household structure, and several types of individuals,
{\it Adv.\ Appl.\ Prob.} {\bf 44} 63--86.

{\sc Ball, F.G.; Sirl, D.J.\ and Trapman, P.} (2009),
Threshold behaviour and final outcome of an epidemic on a random network with household structure,
{\it Adv.\ Appl.\ Prob.} {\bf 41} 765--796.

{\sc Ball, F.G.; Sirl, D.J.\ and Trapman, P.} (2010),
Analysis of a stochastic SIR epidemic on a random network incorporating household structure,
{\it Math.\ Biosci.} {\bf 224(2)} 53--73.

{\sc Ball, F.G.; Sirl, D.J.\ and Trapman, P.} (2012),
Epidemics on random intersection graphs,
Submitted.

{\sc \Barabasi, A.\ and Albert, R.} (1999),
Emergence of scaling in random networks,
{\it Science} {\bf 286} 509--512.

%(?Right reference I assume?) {\sc Bollob{\'a}s, B.} (2001),
%{\it Random Graphs},
%Cambridge University Press, Cambridge.

{\sc Britton T.; Nordvik, M.K.\ and Liljeros, F.} (2007)
Modelling sexually transmitted infections: the effect of partnership activity and number of partners on $R_0$,
{\it Theor.\ Pop.\ Biol.} {\bf 72} 389-399.

{\sc Britton T.; Deijfen, M.; Lindholm, M.\ and \Lageras, A.N.} (2008),
Epidemics on random graphs with tunable clustering,
{\it J.\ Appl.\ Prob.} {\bf 45} 743--756.

{\sc Diekmann, O.\ and Heesterbeek, J.A.P.} (2000),
{\it Mathematical Epidemiology of Infectious Diseases}, Chichester:
John Wiley \& Son.

{\sc Diekmann, O.; de Jong, M.C.M.\ and Metz, J.A.J.} (1998),
A deterministic epidemic model taking account of repeated contacts between the same individuals,
{\it J.\ Appl.\ Prob.} {\bf 35} 448--462.

{\sc \Erdos, P.\ and \Renyi, A.} (1959),
On random graphs. {\it Publicationes Mathematicae} {\bf 6}, 290-297.

{\sc Gleeson, J.P.} (2009),
Bond percolation on a class of clustered random networks,
{\em Phys.\ Rev.\ E} {\bf 80}, 036107.

{\sc Gleeson, J.P.; Melnik, S.\ and Hackett, A.} (2010),
How clustering affects the bond percolation threshold in complex networks,
{\em Phys.\ Rev.\ E} {\bf 81}, 066114.

{\sc van der Hofstad, R.\ and Litvak, N.} (2012),
Degree-degree correlations in random graphs with heavy-tailed degrees,
arXiv:1202.307v3.

{\sc Isham, V., Kaczmarska, J. and Nekovee, M.} (2011),
Spread of information and infection on finite random networks.
{\em Phys.\ Rev.\ E} {\bf 83}, 046128.

{\sc Karrer, B.\ and Newman, M.E.J.} (2010),
Random graphs containing arbitrary distributions of subgraphs,
{\em Phys.\ Rev.\ E} {\bf 82}, 066118.

{\sc Ma, J.; van den Driessche, P.\ and Willeboordse, F.H.} (2012),
Effective degree household network disease model,
{\em J.\ Math.\ Biol.} Published online 18th January 2012. DOI 10.1007/s00285-011-0502-9.

{\sc May, R.M. and Anderson, R.M.} (1987),
Transmission dynamics of HIV infections,
{\em Nature} {\bf 326}, 137--142.

%{\sc Miller, J.C.} (2009),
%Spread of infectious disease through clustered populations,
%{\em J.\ R.\ Soc.\ Interface} {\bf 6} 1121--1134.

{\sc Miller, J.C.} (2009),
Percolation and epidemics in random clustered networks,
{\em Phys.\ Rev.\ E} {\bf 80} 020901(R).

{\sc Mode, C.J.} (1971),
{\em Multitype branching processes. {T}heory and applications}.
Modern Analytic and Computational Methods in Science and Mathematics, {\bf 34}. Elsevier, New York.

{\sc Mollison, D.} (1977),
Spatial contact models for ecological and epidemic spread,
{\em J.\ Roy.\ Stat.\ Soc.\ B} {\bf 39}(3) 283--326.

{\sc Molloy, M. and Reed, B.} (1995),
A critical point for random graphs with a given degree sequence.
{\em Rand.\ Struct.\ Alg.} {\bf 6} 161--179.

{\sc Newman, M.E.J., Strogatz, S.H. and Watts, D.J.} (2001),
Random graphs with arbitrary degree distributions and their applications,
{\it Phys.\ Rev.\ E} {\bf 64}, 026118.

{\sc Newman, M.E.J.} (2002a), Assortative mixing in networks,
{\it Phys.\ Rev.\ Lett.} {\bf 89} 208701.

{\sc Newman, M.E.J.} (2002b),
Spread of epidemic disease on networks,
{\it Phys.\ Rev.\ E} {\bf 66} 016128.

{\sc Newman, M.E.J.} (2003),
The structure and function of complex networks,
{\it SIAM Review} {\bf 45} 167--256.

{\sc Newman, M.E.J.} (2009),
Random graphs with clustering,
{\it Phys.\ Rev.\ Lett.} {\bf 103} 058701.

{\sc Trapman, P.} (2007),
On analytical approaches to epidemics on networks,
{\it Theor.\ Pop.\ Biol.} {\bf 71} 160--173.

{\sc Watts, S.C.\ and Strogatz, S.H.} (1998),
Collective dynamics of `small-world' networks,
{\em Nature} {\bf 393} 440--442.

\end{document}